\documentclass[12pt]{amsart}

\usepackage{
	amsmath,
	amsfonts,
	amssymb,
	amsthm,
	amscd,
	booktabs,
	comment,
	enumitem,
	etoolbox,
	gensymb,    
	mathtools,
	mathdots,
	stmaryrd,
}
\usepackage[usenames,dvipsnames]{xcolor}
\usepackage[all]{xy}

\DeclareMathOperator{\Hom}{Hom}
\DeclareMathOperator{\End}{End}
\newcommand*{\SL}{\mathrm{SL}}
\newcommand*{\bQ}{\mathbb{Q}}
\newcommand*{\bC}{\mathbb{C}}
\newcommand*{\fg}{\mathfrak{g}}
\newcommand*{\fS}{\mathfrak{S}}
\newcommand*{\Dfn}[1]{\emph{\color{blue}#1}}
\newtheorem{thm}{Theorem}[section]
\newtheorem{defn}[thm]{Definition}
\newtheorem{lemma}[thm]{Lemma}
\newtheorem{prop}[thm]{Proposition}

\newtheorem{cor}[thm]{Corollary}

\usepackage{tikz-cd}

\usepackage[T1]{fontenc}
\usepackage[colorlinks=true, linkcolor=blue, citecolor=blue, urlcolor=blue, breaklinks=true]{hyperref}

\DeclareFontFamily{OT1}{pzc}{}
\DeclareFontShape{OT1}{pzc}{m}{it}{<-> s * [1.10] pzcmi7t}{}
\DeclareMathAlphabet{\mathpzc}{OT1}{pzc}{m}{it}

\newcommand\Ccat{\mathpzc{C}}
\newcommand\Gcat{\mathpzc{G}}
\newcommand\Fcat{\mathpzc{F}}
\newcommand\Tcat{\mathpzc{T}}
\newcommand\Dcat{\mathpzc{D}}
\newcommand\kk{\Bbbk}
\newcommand\Mcat{\mathpzc{M}}
\newcommand\Pcat{\mathpzc{P}}
 
\newcommand\bB{\mathcal{B}}
\newcommand\Rot{\mathrm{Rot}}

\newcommand\Tr{\mathrm{Tr}}

\newcommand{\braid}{\boldsymbol\beta}

\newcommand{\co}{\mathrm{co}}
\newcommand{\ev}{\mathrm{ev}}

\newtheorem*{assume}{Assumption}

\tikzstyle{repV}=[red, ultra thick]
\tikzstyle{repW}=[black, ultra thick, dashed]
\tikzstyle{bgd}=[blue, fill=blue!20]
\newcommand\posdot[1]{
	\filldraw[black] (#1) -- +(0.09192,0.09192) arc(45:135:0.13) -- cycle;
	\filldraw[black] (#1) -- +(0.09192,-0.09192) arc(-45:-135:0.13) -- cycle;
}
\newcommand\negdot[1]{
	\filldraw[black] (#1) -- +(0.09192,0.09192) arc(45:-45:0.13) -- cycle;
	\filldraw[black] (#1) -- +(-0.09192,0.09192) arc(135:225:0.13) -- cycle;
}

\begin{document}

\title{A universal plane of diagrammatic categories}

\author{Bruce W.\ Westbury}
\email{brucewestbury@gmail.com}

\begin{abstract}
	We introduce quantum skein relations for a family of ribbon categories parametrised by the projective plane over $\bQ$.
	There are thirteen points for which the ribbon category admits a ribbon functor to a category of invariant tensors
	for a quantised  enveloping algebra. These thirteen points lie on three projective lines. One line gives the first
	row of the Freudenthal magic square and another gives the fourth row which is the exceptional series.
\end{abstract}

\subjclass[2020]{18M15, 18M30, 17B37}

\keywords{monoidal category, quantized enveloping algebra}

\maketitle
\thispagestyle{empty}

\tableofcontents

\section{Introduction}
One of the main objectives in tensor invariant theory is to describe the \emph{category of invariant tensors on $V$}; this is the full pivotal subcategory on $V$ and its dual.
Originally $V$ was taken to be a rational representation of a reductive algebraic group.
This problem was first raised in \cite{Littlewood1944} and the influential
book \cite{Weyl1947}. More recently this problem has been studied for $V$ a finite
dimensional representation of a quantised enveloping algebra. The first example to be
studied was $U_q(\SL(2))$, using the diagrams of \cite{Levinson1962}. Subsequently the fundamental representations of rank two
simple Lie algebras were studied in \cite{Kuperberg1997}; the spin representation of $B_3$ was studied in \cite{Westbury2008}; the fundamental representations
of $U_q(\SL(n))$ were studied in \cite{Cautis2014}; and the fundamental representations of $U_q(\mathrm{Sp}(2n))$ are studied in \cite{Bodish2021}.

We give a uniform diagrammatic approach to several categories of invariant tensors. The motivation in the short term  is simply to extend the list of examples.
The motivation in the long term is to give some structure to skein theories with some organising principles. This is similar to the motivation for \cite{Morrison2016}
although the technical details are different. The next step in developing skein theories would be to apply the same methods to the Vogel plane,
\cite{Vogel2011}. The ultimate aim would be develop skein theories for all representations $V$ such that $V$ is self-dual and $V\otimes V$ is multiplicity free.

The representations we consider are highest weight representations of a semisimple Lie algebra, $L(\lambda)$, such that the symmetric and exterior square decompose as
\begin{equation*}
	S^2 L(\lambda) \cong L(0) \oplus L(2\lambda) \oplus L(\mu)
	\qquad
	\bigwedge\nolimits^2 L(\lambda) \cong L(\theta) \oplus L(\nu)
\end{equation*}
for some dominant weights $\mu$ and $\nu$.
Here $0$ is the zero weight so $L(0)$ is the trivial representation and $\theta$ is the highest root so $L(\theta)$ is the adjoint representation.
There are eight examples for representations of classical groups shown in Figure~\ref{fig:classical} and 
six examples for representations of exceptional groups shown in Figure~\ref{fig:exceptional}.

\begin{figure}[ht]
\begin{tabular}{c|cccccccc}
	$G$ & $A_1$ & $A_2$ & $A_7$ & $B_4$ & $C_3$ & $C_4$ & $D_4$ & $D_8$ \\
	$\lambda$ & $4\omega_1$ & $\omega_1+\omega_2$ & $\omega_4$ & $\omega_4$ & $\omega_2$ & $\omega_4$ & $\omega_2$ & $\omega_8$ \\
	$\dim$ & 5 & 8 & 16 & 70 & 14 & 42 & 28 & 128
\end{tabular}
\caption{Classical cases}\label{fig:classical}
\end{figure}

\begin{figure}[ht]
\begin{tabular}{c|cccccc}
	 $G$ & $G_2$ & $F_4$ & $F_4$ & $E_6$ & $E_7$ & $E_8$ \\
	 $\lambda$ & $\omega_2$ & $\omega_1$ & $\omega_4$ & $\omega_2$ & $\omega_1$ & $\omega_8$ \\
	$\dim$ & 14 & 26 & 52 & 78 & 133 & 248
\end{tabular}
\caption{Exceptional cases}\label{fig:exceptional}
\end{figure}

Our main result is an extension of \cite{Savage2022} and we also follow this paper. The main result is to give the structure constants of $\End(L(\lambda)\otimes L(\lambda))$.
This is a five dimensional commutative algebra and the structure constants are with respect to a basis of diagrams. It is unlikely that these relations are complete
for any of the examples in Figures~\ref{fig:classical} and \ref{fig:exceptional}.

In order to present an overview of this paper I will start with Section~\ref{sec:classical}
and work backwards. In Section~\ref{sec:classical} we give a finite presentation of a pivotal symmetric category. In a first reading, this category 
can be taken as linear over $\bQ^{(0)}(m,n,p)$, the ring of homogeneous rational functions of degree 0 in the three indeterminates, $m,n,p$.
However this does not allow specialisations. Let $R_{cl}$ be the finitely generated subring of $\bQ^{(0)}(m,n,p)$ generated by the coefficients of the skein relations.
A homomorphism $R_{cl}\to\bQ$ is determined by giving the values of $m,n,p$
and almost all values determine a homomorphism. Given a ring homomorphism $R_{cl}\to\bQ$ we can specialise to give a pivotal symmetric $\bQ$-linear
category. The intuition is that we have a family of pivotal symmetric $\bQ$-linear categories parametrised by a Zariski open subset of the projective plane over $\bQ$.
The property of the classical category is that, for each example in Figures~\ref{fig:classical} and \ref{fig:exceptional}
there is a specialisation $R_{cl}\to\bQ$ with a pivotal symmetric $\bQ$-linear functor to the category of invariant tensors for the Lie algebra.

There are several examples of families of pivotal symmetric $\bQ$-linear categories parametrised by a Zariski open subset of the projective line over $\bQ$.
The original example is the Brauer category, \cite{Brauer1937}; and there is also the oriented version of this category.
Further examples are given in the pioneering book, \cite[Chapters~16--20]{Cvitanovic2008}; these are associated to rows of the Freudenthal
magic square, \cite{Freudenthal1964}. Our category specialises to the first and fourth rows. The only example of a family of pivotal symmetric $\bQ$-linear categories
parametrised by a Zariski open subset of the projective plane over $\bQ$ is the Vogel plane constructed in \cite{Vogel2011}.
A different approach is the theory of interpolation categories presented in \cite{Deligne2007}.

The relations in Section~\ref{sec:classical} are not derived directly but are obtained by taking the classical limit of quantum relations.
This runs counter to the traditional approach of constructing the classical version first and then quantising. In Section~\ref{algebra}
we give a finite presentation of a ribbon category.
In Section~\ref{scalars} we construct a domain, $D(t,u,v)$,
whose field of fractions is the field of rational functions in $t,u,v$; together with a homomorphism $D(t,u,v)\to\bQ^{(0)}(m,n,p)$.
In a first reading, the coefficients of the skein relations in Section~\ref{algebra} are elements of $D(t,u,v)$.
The relations in Section~\ref{sec:classical} are obtained by applying ring homomorphism $D(t,u,v)\to\bQ^{(0)}(m,n,p)$ to the skein relations in Section~\ref{algebra}.

However this does not allow specialisations. Let $R_{qu}$ be the finitely generated subring of $D(t,u,v)$ generated by the coefficients of the skein relations.
Then the homomorphism $D(t,u,v)\to\bQ^{(0)}(m,n,p)$ restricts to $R_{qu}\to R_{cl}$.
The property of the quantum category is that, for each example in Figures~\ref{fig:classical} and \ref{fig:exceptional}
there is a specialisation $R_{qu}\to\bQ(q)$ with a $\bQ(q)$-linear ribbon functor to the category of invariant tensors for the quantised enveloping algebra.
The construction of these functors
is based on the formula for the action of the ribbon element on $L_q(\lambda)$ and on the formula for the quantum dimension of $L_q(\lambda)$
which is a quantum analogue of the Weyl dimension formula. In order to show the functor satisfies the relations, we show in Section~\ref{sec:relations}, that
the relations are determined by this information.

The only families which have quantum versions are the Brauer category whose quantisation is the BMW category, \cite{Birman1989},
and the oriented version \cite{Brundan2017}. These families are parametrised by a projective line.
The category we define is the first example of a quantum analogue of a family of pivotal symmetric categories parametrised
by a plane. Our category specialises to a quantum analogue of the pivotal symmetric categories associated to the first and fourth rows of the
Freudenthal magic square. These are both new results.

We now discuss some features of the family of pivotal symmetric categories in Section~\ref{sec:classical}. These features were unexpected although similar
features appear in the Vogel plane, \cite{Vogel2011}. The examples in Figures~\ref{fig:classical} and \ref{fig:exceptional} lie on three lines in the plane.
One line corresponds to the first row of the Freudenthal magic square, a second line corresponds to the fourth row of the Freudenthal magic square
and the third line is the series of exceptional symmetric spaces first observed in \cite[\S~8.4]{Westbury2015}. This family also admits an involution.
There are also three degenerate lines. One gives the Birman-Wenzl category (after semisimplification), another has all points corresponding to $\SL(2)$
and the interpretation of the third degenerate line is not clear.

\subsubsection*{SageMath computations}

Throughout the paper, many direct computations are performed using the open-source mathematics software system SageMath, \cite{sagemath}.

\section{Ribbon categories}

\subsection{Monoidal categories}
There is a host of definitions of a monoidal category in \cite{Leinster2004}. Here we give the most commonly used definition from \cite{MacLane1963}.

\begin{defn} The data for a \Dfn{monoidal category} (aka tensor category) consists of;
a category $C$, a functor $\otimes\colon C\times C\to C$, an object $I$ of $C$
together with natural isomorphisms
$\alpha_{X,Y,Z}\colon (X\otimes Y)\otimes Z\to X\otimes (Y\otimes Z)$,
natural isomorphisms $\lambda_X\colon X\otimes I\to X$ and natural isomorphisms
$\rho_X\colon I\otimes X\to X$. 
	
The natural transformations are required to satisfy the pentagon identity  and the triangle identity.
\end{defn}

A monoidal category is \Dfn{strict} if the constraints are identity natural transformations.

If the category $C$ is essentially small then the set of isomorphism classes of objects is a monoid.

\begin{defn} A \Dfn{strong monoidal functor} consists of a functor $F\colon C\to D$,
	natural isomorphisms $\phi_{X,Y} \colon F(X)\otimes F(Y) \to F(X\otimes Y))$
	and an isomorphism $I\to F(I)$. The natural isomorphisms satisfy the constraint that
	the following diagrams commute
	\begin{equation*}
		\begin{tikzcd}
			(F(X)\otimes F(Y))\otimes F(Z) \arrow[r] \arrow[d, swap, "\phi_{X,Y}\otimes 1_Z"] & F(X) \otimes (F(Y)\otimes F(Z)) \arrow[d, "1_X\otimes\phi_{Y,Z}"] \\
			F(X\otimes Y)\otimes F(Z) \arrow[d, swap, "\phi_{X\otimes Y,Z}"] & F(X) \otimes F(Y\otimes Z) \arrow[d, "\phi_{X,Y\otimes Z}"] \\
			F((X\otimes Y)\otimes Z) \arrow[r] & F(X\otimes (Y\otimes Z))
		\end{tikzcd}
	\end{equation*}
	
	\begin{equation*}
		\begin{tikzcd}
			F(X)\otimes I \arrow[r] \arrow[d] & F(X) \arrow[d] \\
			F(X)\otimes F(I) \arrow[r] & F(X\otimes I)
		\end{tikzcd}
		\qquad
		\begin{tikzcd}
			I \otimes F(X) \arrow[r] \arrow[d] & F(X) \arrow[d] \\
			F(I)\otimes F(X) \arrow[r] & F(I\otimes X)
		\end{tikzcd}
	\end{equation*}
\end{defn}

\begin{defn} A \Dfn{monoidal transformation} $F\to G$ is a natural transformation
	$\eta\colon F\to G$ such that the following diagrams commute
	\begin{equation*}
		\begin{tikzcd}
			F(X)\otimes F(Y) \arrow[r] \arrow[d, swap, "\eta_X\otimes\eta_Y"] &
			F(X\otimes Y) \arrow[d, "\eta_{X\otimes Y}"] \\
			G(X)\otimes G(Y) \arrow[r] & G(X\otimes Y)
		\end{tikzcd}
		\qquad
		\begin{tikzcd}[row sep=tiny]
			& F(I) \arrow[dd, "\eta_I"] \\
			I \arrow[ur] \arrow[dr] & \\
			& G(I)
		\end{tikzcd}
	\end{equation*}	
\end{defn}

\subsection{Braided categories} In this section we discuss braided and balanced
categories. Braided categories were introduced in \cite{Joyal1993}.
Braided and balanced categories are studied in \cite{Joyal1991a}.

\begin{defn} A \Dfn{braided category} is a monoidal category with natural isomorphisms
	$\braid_{X,Y}\colon X\otimes Y \to Y\otimes X$ such that the hexagon identities are satisfied.
\end{defn}

\begin{defn} A \Dfn{braided functor} is a monoidal functor such that the following diagram
	commutes
	\begin{equation*}
		\begin{tikzcd}
			F(X)\otimes F(Y) \arrow[r] \arrow[d, swap, "\braid_{F(X),F(Y)}"] &
			F(X\otimes Y) \arrow[d, "F(\braid_{X,Y})"] \\
			F(Y) \otimes F(X) \arrow[r]& F(Y\otimes X)
		\end{tikzcd}
	\end{equation*}
\end{defn}

\begin{lemma} Given a braided category, $\Ccat$, $(1,\braid^2,1)$ is a monoidal functor $\Ccat\to \Ccat$.
\end{lemma}
Taken from \cite[\S~1]{Shum1994}. 

\begin{defn} A \Dfn{twist} is a monoidal isomorphism 
	\begin{equation*}
		\theta \colon (1,1,1)\to (1,\braid^2,1)
	\end{equation*}
\end{defn}
This means that
\begin{equation*}
	\theta_{X\otimes Y} = \braid_{Y,X}\braid_{X,Y} (\theta_X\otimes\theta_Y)
\end{equation*}

\begin{defn} A \Dfn{balanced category} is a braided category with a twist.
\end{defn}

\begin{defn} A \Dfn{balanced functor} is a braided functor such that
	$F(\theta_X) = \theta_{F(X)}$ for all objects $X$.
\end{defn}

\subsection{Pivotal categories}
In this section we discuss duals, autonomous categories and pivotal categories.

Let $X$ be an object of a monoidal category $C$. A \Dfn{right dual} of $X$
is an object $X^\ast$ of $C$ together with morphisms $\co_X\colon I\to X^\ast\otimes X$
and $\ev_X\colon X\otimes X^\ast\to I$ such that the composites in \eqref{zigzag1} and \eqref{zigzag2}
are the identity.
\begin{equation}\label{zigzag1}
	X\to X\otimes I\to X\otimes (X^\ast\otimes X)	\to (X\otimes X^\ast)\otimes X
	\to I\otimes X \to X
\end{equation}
\begin{equation}\label{zigzag2}
	X^\ast \to I\otimes X^\ast \to (X^\ast\otimes X) \otimes X^\ast\to 
	X^\ast (X\otimes X^\ast) \to X^\ast\otimes I \to X^\ast	
\end{equation}
This will be denoted by  $(\co_X,\ev_X)\colon X\dashv X^\ast$.

Suppose $F$ is a strong monoidal functor and $(\co_X,\ev_X)\colon X\dashv X^\ast$
is an adjunction. Then we have an adjunction for $F(X^\ast)$ given by the composites
in \eqref{Fev} and \eqref{Fco}.
\begin{equation}\label{Fev}
	I\to F(I)\to F(X^\ast\otimes X) \to F(X^\ast)\otimes F(X)	
\end{equation}
\begin{equation}\label{Fco}
	F(X)\otimes F(X^\ast) \to F(X\otimes X^\ast) \to F(I) \to I	
\end{equation}

\begin{defn} An \Dfn{autonomous category} is a monoidal category together
	with an adjunction $(\co_X,\ev_X)\colon X\dashv X^\ast$ for all objects $X$.
\end{defn}

There is no need for a definition of autonomous functor or autonomous transformation
since right duals are unique up to isomorphism.

Let $f\colon X\to Y$ be a morphism in a monoidal category. Suppose $(X^\ast,\co_X,\ev_X)$
is a right dual for $X$ and $(Y^\ast,\co_Y,\ev_Y)$ is a right dual for $Y$. 
Then we define a map $f^\ast\colon Y^\ast\to X^\ast$ to be the composite
in \eqref{dualmap}.
\begin{equation}\label{dualmap}
	Y^\ast \to I\otimes Y^\ast\to (X^\ast\otimes X)\otimes Y^\ast
	\to (X^\ast\otimes Y) \otimes Y^\ast \to X^\ast \otimes (Y\otimes Y^\ast)
	\to X^\ast\otimes I \to X^\ast
\end{equation}

If $C$ is autonomous then $f\mapsto f^\ast$ is a functor $\Ccat\to \Ccat^{\mathrm{op}}$.
Then we should check the identity $(f\otimes g)^\ast = g^\ast \otimes f^\ast$.
Then $\ast$ is a strong monoidal functor.

The following is proved in \cite[Proposition~5.2.3]{rivano2006};
see also \cite[Proposition~7.1]{Joyal1993}.
\begin{lemma} Suppose $\tau\colon F\to G$ is a monoidal transformation between
	monoidal functors $F,G\colon C\to D$. If $X$ has a	right dual $X^\ast$ in $C$
	then  the following diagram commutes
	\begin{equation*}
		\begin{tikzcd}
			F(X^\ast) \arrow[r, "\tau_{X^\ast}"] \arrow[d] & G(X^\ast) \arrow[d] \\
			(F(X))^\ast \arrow[r, swap, "(\tau_X)^\ast"] & (G(X))^\ast
		\end{tikzcd}
	\end{equation*}
\end{lemma}

The definition of a sovereign category is \cite[Definition~4.1]{freyd1992}.
\begin{defn} A \Dfn{pivotal category} (aka sovereign category) is an autonomous category together
	with natural isomorphisms $i_X\colon X\to X^{\ast\ast}$.
\end{defn}

The definition of a ribbon category is given in \cite[\S~6,7]{Joyal1993} and  \cite{Shum1994}.
The motivation for the definition is given in \cite[A.2]{Henriques2016}.

\begin{defn} A \Dfn{ribbon category} (aka tortile category) is both balanced
	and pivotal. The twist and pivotal structures satisfy a condition which implies
	that each determines the other.
\end{defn}

\begin{defn} A \Dfn{ribbon category} is a pivotal, balanced category
	whose twist $\theta$ satisfies the additional condition
	\begin{equation*}
		\theta_{V^*} = (\theta_V)^*
	\end{equation*}
\end{defn}

\begin{defn} A \Dfn{rigid symmetric} category is a ribbon category in which
	the braiding satisfies $\braid_{X,Y}\braid_{Y,X}=1_{X\otimes Y}$, for all $X,Y$
	and the twist satisfies $\theta_X=1_X$ for all $X$.
\end{defn}

\subsection{Free monoidal categories}

\begin{defn}
Let $\Ccat$ be a monoidal category. The category $V(\Ccat)$ has objects consisting
of a pair of self duals objects, $V$ and $W$, and a morphism $T\colon V\otimes V\to W$.
The data that makes $V$ and $W$  self-dual are adjunctions $(\co_V,\ev_V)\colon V\dashv V$
and $(\co_W,\ev_W)\colon W\dashv W$. A morphism $(V,W,T)\to (V',W',T')$ consists of a
pair of morphisms $f_V\colon V\to V'$ and $f_W\colon W\to W'$. Then $f_V$ gives
the adjunction $(f(\co_V),f(\ev_V))\colon V'\dashv V'$ and similarly for $W$.
Then we also require that $(f_V\otimes f_V)\circ T' = T\circ f_W$.
\end{defn}
We can draw the string diagrams. This plays the role of a generating set.

The map $\Ccat\mapsto V(\Ccat)$ is a functor from the bicategory consisting of monoidal categories,
monoidal functors and monoidal transformations
to the bicategory consisting of categories,
functors and natural transformations.

Suppose $\Ccat$ and $\Ccat'$ are monoidal categories and $N$ is an object of
$V(\Ccat)$. Then there is an induced functor $\Hom(\Ccat,\Ccat')\to V(\Ccat')$.
The map of objects maps a monoidal functor $F\colon \Ccat\to\Ccat'$ to $F(N)$.
The map of morphisms maps a monoidal transformation $\eta\colon F\to G$
to the pair $(\eta_V,\eta_W)$.

The following definition is an instance of \cite[\S~4~Definition~1.5]{Joyal1991a}.
\begin{defn}  Suppose $(\Mcat,N)$ is a pair consisting of a monoidal category
	$\Mcat$ and an object $N$ of $V(\Mcat)$. A \Dfn{free monoidal category} is a pair
	$(\Mcat,N)$ such that
	the induced functor $\Hom(\Mcat,\Ccat)\to V(\Ccat)$ is an equivalence of categories, for all  monoidal categories $\Ccat$.
\end{defn}
If $(\Mcat,N)$ is a free monoidal category then the induced functor $\Hom(\Mcat,\Ccat)\to V(\Ccat)$ is an isomorphism of categories, for all strict monoidal categories $\Ccat$.

This definition has a natural modification for monoidal categories with structure.
In section \ref{sec:Tdef} we give a free rigid symmetric category in Theorem~\ref{thm:freerigid},
a free pivotal category in Theorem~\ref{thm:freepivotal} and a free ribbon category in Theorem~\ref{thm:freeribbon}.

\section{Diagram categories}\label{sec:Tdef}
In this section we construct a free rigid symmetric monoidal category
and a free ribbon category. Both of these categories are constructed as diagram categories.

\subsection{Free symmetric category}\label{sec:Gcat}
First we introduce the classical diagram category. This is a strict rigid symmetric monoidal category. This is a cobordism category of labelled trivalent graphs.

\begin{defn} A \Dfn{trivalent graph}, $\Gamma$, is a graph whose vertices all have degree one or three,
\begin{itemize}
	\item all edges are labelled $V$ or $W$
	\item each trivalent vertex has one edge labelled $W$ and two labelled $V$
\end{itemize} 
The boundary of a trivalent
graph $\Gamma$ is $\partial\Gamma$, the set of univalent vertices of $\Gamma$.
\end{defn}
When we draw a trivalent graph we will colour the edges instead of labelling them.

\begin{defn} An object of the category $\Gcat$ is a finite ordered set $X$ together with a function $X\to \{V,W\}$, called the labelling.
	A morphism $X\to Y$ consists of a trivalent graph $\Gamma$
together with an isomorphism of labelled sets $\partial\Gamma\cong X\amalg Y$. Composition is given by glueing,
the identity object of a labelled set $X$ is $X\times I$.

This category is a strict rigid symmetric category. The tensor product is given by disjoint union.
Every object is self-dual and the cup and caps are shown in \eqref{gen:pivotal}. The symmetric structure is given by reordering the finite ordered sets.
\end{defn}

Suppose $\Ccat$ and $\Ccat'$ are rigid symmetric categories and $N$ is an object of
$V(\Ccat)$. Then there is an induced functor $\Hom(\Ccat,\Ccat')\to V(\Ccat')$.
The map of objects maps a rigid symmetric functor $F\colon \Ccat\to\Ccat'$ to $F(N)$.
The map of morphisms maps a rigid symmetric transformation $\eta\colon F\to G$
to the pair $(\eta_V,\eta_W)$.

\begin{defn}  Suppose $(\Mcat,N)$ is a pair consisting of a rigid symmetric category
	$\Mcat$ and an object $N$ of $V(\Mcat)$. A \Dfn{free rigid symmetric category} is a pair
	$(\Mcat,N)$ such that
	the induced functor $\Hom(\Mcat,\Ccat)\to V(\Ccat)$ is an equivalence of categories, for all  rigid symmetric categories $\Ccat$.
\end{defn}

There is an obvious object $N$ of $V(\Gcat)$.
\begin{thm}\label{thm:freerigid}
	The pair $(\Gcat,N)$ is a free rigid symmetric category.
\end{thm}

\subsection{Free pivotal category}
In this section we give a topological construction of a free pivotal category.

\begin{figure}[ht]
\begin{equation*}
\begin{tikzpicture}
\filldraw[bgd] (1,0) -- (-1,0) arc(180:360:1);
\draw[black,thick] (1,0) -- (-1,0);
\draw[repV] (0,0) -- (0,-1);
\end{tikzpicture}
\qquad
\begin{tikzpicture}
	\filldraw[bgd] (1,0) -- (-1,0) arc(180:360:1);
	\draw[black,thick] (1,0) -- (-1,0);
	\draw[repW] (0,0) -- (0,-1);
\end{tikzpicture}
\end{equation*}
\begin{equation*}
\begin{tikzpicture}
	\filldraw[bgd] (-1,0) -- (1,0) arc(0:180:1);
	\draw[black,thick] (1,0) -- (-1,0);
	\draw[repV] (0,0) -- (0,1);
\end{tikzpicture}
\qquad
\begin{tikzpicture}
	\filldraw[bgd] (-1,0) -- (1,0) arc(0:180:1);
	\draw[black,thick] (1,0) -- (-1,0);
	\draw[repW] (0,0) -- (0,1);
\end{tikzpicture}
\end{equation*}
\caption{Boundary points}\label{fig:edge}
\end{figure}
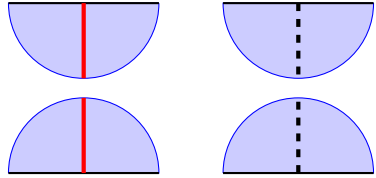

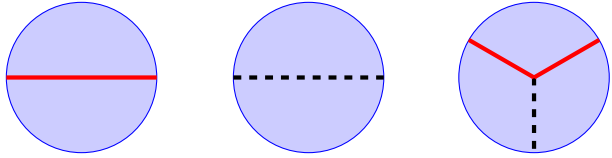
\begin{figure}[ht]
	\begin{tikzpicture}[baseline=0]
	\draw [bgd] (0,0) circle [radius=1.0];
	\draw [repV] (-1,0) -- (1,0);
\end{tikzpicture}
\qquad
	\begin{tikzpicture}[baseline=0]
	\draw [bgd] (0,0) circle [radius=1.0];
	\draw [repW] (-1,0) -- (1,0);
\end{tikzpicture}
\qquad
	\begin{tikzpicture}[baseline=0]
	\draw [bgd] (0,0) circle [radius=1.0];
	\draw [repW] (0,0) -- (270:1);
	\draw [repV] (0,0) -- (30:1);
	\draw [repV] (0,0) -- (150:1);
\end{tikzpicture}
	\caption{Interior points}\label{fig:interior}
\end{figure}

A \Dfn{rectangle} will mean a rectangle embedded in the Euclidean plane with horizontal
and vertical edges.
\begin{defn} A \Dfn{planar diagram} is a closed subset $L$ of a rectangle such that
every point in $L$ which is in the boundary of the rectangle has an open neighbourhood
equivalent to one of the open neighbourhoods in Figure~\ref{fig:edge} and every point of
$L$ in the interior of the rectangle has an open neighbourhood equivalent
to one of the open neighbourhoods in Figure~\ref{fig:interior}.
\end{defn}

\begin{defn}\label{defn:pivotal}
	 The set of objects of $\Pcat$ is the free monoid on an alphabet with two
	letters, so an element of this monoid is a word in this alphabet. A morphism of the category $\Pcat$ is an equivalence class of planar diagrams. The equivalence relation is ambient
	isotopy preserving the edges of the rectangle as sets. The domain of a planar diagram
	is the set of intersection points of $L$ with the top edge, read from left to right,
	and the codomain of a planar diagram
	is the set of intersection points of $L$ with the bottom edge, again, read from left to right.
	Composition is given, as usual, by stacking rectangles.
	
	This is a strict autonomous category. The tensor product is juxtaposition of rectangles.
	The morphisms in the adjunctions defining duals have no trivalent vertices.
\end{defn}
The functor $\ast\colon \Pcat\to \Pcat^{\mathrm{op}}$ rotates a diagram through a half revolution.
Applying this twice rotates a diagram through one revolution which does not change the diagram
so the category $\Pcat$ is strict pivotal.

Suppose $\Ccat$ and $\Ccat'$ are pivotal categories and $N$ is an object of
$V(\Ccat)$. Then there is an induced functor $\Hom(\Ccat,\Ccat')\to V(\Ccat')$.
The map of objects maps a pivotal functor $F\colon \Ccat\to\Ccat'$ to $F(N)$.
The map of morphisms maps a pivotal transformation $\eta\colon F\to G$
to the pair $(\eta_V,\eta_W)$.

\begin{defn} A \Dfn{free pivotal category} is a pair
	$(\Mcat,N)$ consisting of a pivotal category
	$\Mcat$ and an object $N$ of $V(\Mcat)$ such that
	the induced functor $\Hom(\Mcat,\Ccat)\to V(\Ccat)$ is an equivalence of categories, for all  pivotal categories $\Ccat$.
\end{defn}

There is an obvious object $N$ of $V(\Pcat)$.
\begin{thm}\label{thm:freepivotal}
	The pair $(\Pcat,N)$ is a free pivotal category.
\end{thm}

\subsection{Free ribbon category}
The free ribbon category is constructed
topologically in \cite{Reshetikhin1990} using ribbon graphs.
A ribbon graph is given by
taking the neighbourhood of a graph embedded in an oriented surface and then embedding
this neighbourhood in three dimensional box. The diagrams we draw are planar and
represent a projection of the box to a rectangle with the ribbon graph projecting
to a neighbourhood of the planar graph.

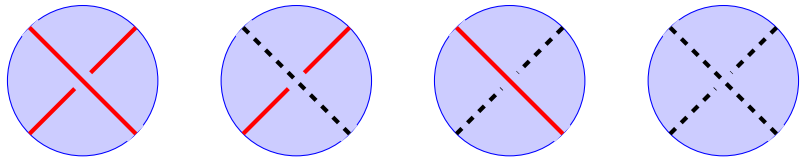
\begin{figure}[ht]
\begin{equation*}
\begin{tikzpicture}[baseline=0]
	\draw [bgd] (0,0) circle [radius=1.0];
	\draw [repV] (45:1) -- (225:1);
	\draw [blue!20, line width=3mm] (135:1) -- (315:1);
	\draw [repV] (135:1) -- (315:1);
\end{tikzpicture}
\qquad
\begin{tikzpicture}[baseline=0]
	\draw [bgd] (0,0) circle [radius=1.0];
	\draw [repV] (45:1) -- (225:1);
	\draw [blue!20, line width=3mm] (135:1) -- (315:1);
	\draw [repW] (135:1) -- (315:1);
\end{tikzpicture}
\qquad
\begin{tikzpicture}[baseline=0]
	\draw [bgd] (0,0) circle [radius=1.0];
	\draw [repW] (45:1) -- (225:1);
	\draw [blue!20, line width=3mm] (135:1) -- (315:1);
	\draw [repV] (135:1) -- (315:1);
\end{tikzpicture}
\qquad
\begin{tikzpicture}[baseline=0]
	\draw [bgd] (0,0) circle [radius=1.0];
	\draw [repW] (45:1) -- (225:1);
	\draw [blue!20, line width=3mm] (135:1) -- (315:1);
	\draw [repW] (135:1) -- (315:1);
\end{tikzpicture}
\end{equation*}
\caption{Crossings}\label{fig:crossings}
\end{figure}

We construct a strict ribbon category, $\Fcat$. First we extend the definition of a
planar diagram to allow the additional neighbourhoods shown in Figure~\ref{fig:crossings}
for points of $L$ in the interior of the rectangle.

\begin{figure}[ht]
	\begin{tikzpicture}[baseline=0]
		\draw[bgd] (-0.25,-2.5) rectangle (1.75,1.5);
		\draw[repV] (0.25,1.5) -- (0.25,1);
		\draw[repV] (0.25,1) to[out=270,in=90] (0.75,0);
		\draw [blue!20, line width =3mm](0.75,1) to[out=270,in=90] (0.25,0);
		\draw[repV] (0.75,1) to[out=270,in=90] (0.25,0);
		\draw[repV] (0.75,0) arc(180:360:0.25) -- (1.25,1) arc(0:180:0.25);
		\draw[repV] (0.25,0) -- (0.25,-1);
		\draw[repV] (0.75,-1) to[out=270,in=90] (0.25,-2);
		\draw [blue!20, line width =3mm](0.25,-1) to[out=270,in=90] (0.75,-2);
		\draw[repV] (0.25,-1) to[out=270,in=90] (0.75,-2);
		\draw[repV] (0.75,-2) arc(180:360:0.25) -- (1.25,-1) arc(0:180:0.25);
		\draw[repV] (0.25,-2) -- (0.25,-2.5);
	\end{tikzpicture}
	\quad=\quad
	\begin{tikzpicture}[baseline=0]
		\draw[bgd] (-0.25,-2.5) rectangle (0.75,1.5);
		\draw[repV] (0.25,1.5) -- (0.25,-2.5);
	\end{tikzpicture}
	\qquad
	\begin{tikzpicture}[baseline=0]
		\draw[bgd] (-0.25,-2.5) rectangle (1.75,1.5);
		\draw[repW] (0.25,1.5) -- (0.25,1);
		\draw[repW] (0.25,1) to[out=270,in=90] (0.75,0);
		\draw [blue!20, line width =3mm](0.75,1) to[out=270,in=90] (0.25,0);
		\draw[repW] (0.75,1) to[out=270,in=90] (0.25,0);
		\draw[repW] (0.75,0) arc(180:360:0.25) -- (1.25,1) arc(0:180:0.25);
		\draw[repW] (0.25,0) -- (0.25,-1);
		\draw[repW] (0.75,-1) to[out=270,in=90] (0.25,-2);
		\draw [blue!20, line width =3mm](0.25,-1) to[out=270,in=90] (0.75,-2);
		\draw[repW] (0.25,-1) to[out=270,in=90] (0.75,-2);
		\draw[repW] (0.75,-2) arc(180:360:0.25) -- (1.25,-1) arc(0:180:0.25);
		\draw[repW] (0.25,-2) -- (0.25,-2.5);
	\end{tikzpicture}
	\quad=\quad
	\begin{tikzpicture}[baseline=0]
		\draw[bgd] (-0.25,-2.5) rectangle (0.75,1.5);
		\draw[repW] (0.25,1.5) -- (0.25,-2.5);
	\end{tikzpicture}
	\caption{Ribbon relation}\label{fig:ribbon}
\end{figure}
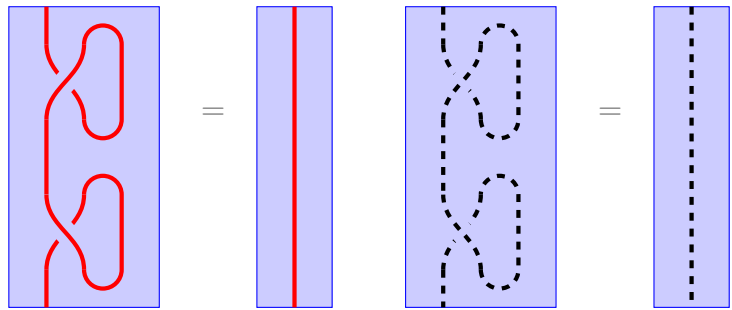

 Define two extended diagrams
to be equivalent under ambient isotopy and the moves $\Omega.2$, $\Omega.3$, $\Omega.3V^-$, $\Omega.3V^+$ in \cite[Definition~2.7]{yetter1989} with all colours for the edges.
The only case of the moves $\Omega.V^\pm$ is where the vertex is the trivalent vertex.
Finally we add the two additional moves shown in Figure~\ref{fig:ribbon}.

Suppose $\Ccat$ and $\Ccat'$ are ribbon categories and $N$ is an object of
$V(\Ccat)$. Then there is an induced functor $\Hom(\Ccat,\Ccat')\to V(\Ccat')$.
The map of objects maps a ribbon functor $F\colon \Ccat\to\Ccat'$ to $F(N)$.
The map of morphisms maps a ribbon transformation $\eta\colon F\to G$
to the pair $(\eta_V,\eta_W)$.

\begin{defn} A \Dfn{free ribbon category} is a pair
	$(\Mcat,N)$ consisting of a ribbon category
	$\Mcat$ and an object $N$ of $V(\Mcat)$ such that
	the induced functor $\Hom(\Mcat,\Ccat)\to V(\Ccat)$ is an equivalence of categories, for all  ribbon categories $\Ccat$.
\end{defn}

There is an obvious object $N$ of $V(\Fcat)$.
\begin{thm}\label{thm:freeribbon}
	The pair $(\Fcat,N)$ is a free ribbon category.
\end{thm}

\subsection{Initial relations}
In this section we assume:
\begin{assume}
$\End(V)=\kk$,  $V\ne I$, $\End(W)=\kk$, and $\Hom(V\otimes V,W)=\kk$.
\end{assume}

The assumption $\End(W)=\kk$ implies:
\begin{equation}\label{tadpole}
	\begin{tikzpicture}[baseline=0]
		\draw [bgd] (0,0) circle [radius=1.0];
		\draw [repW] (0,0.5) -- (0,1);
		\draw [repV] (0,0.5) to [out=210,in=180] (0,-0.5) to [out=0,in=330] (0,0.5);
	\end{tikzpicture}
	\parbox{2cm}{\quad $= 0$}
\end{equation}

Let $\kk$ be a commutative ring. Then we construct a $\kk$-linear ribbon category
$$\Tcat=\Tcat(\alpha,\beta,\gamma,\delta_V,\delta_W,\phi,\psi)$$
where all parameters are in $\kk$ and invertible.

This Proposition assumes we have summands.
\begin{prop}\label{prop:Reigen} Let $X$, $Y$, $Z$ be objects whose endomorphism algebra is $\kk$.
	Assume that $Z$ is a summand of $X\otimes Y$ with multiplicity one..
	Then $\braid_{X,Y}\braid_{Y,X}$ acts on the summand $Z$ by the scalar 
	$\theta_X\theta_Y\theta_Z^{-1}$.
\end{prop}
This is also \cite[(2.21)~Proposition~(2)]{Leduc1997}.
The quantum group version is \cite[(1.37)]{Reshetikhin1987}.

\subsubsection{Initial relations}\label{initial}
We define $\delta_V$, $\delta_W$ to be the quantum dimensions of $V,W$.

\begin{equation*}
	\begin{tikzpicture}[baseline=0]
		\draw [bgd] (0,0) circle [radius=1.0];
		\draw [repV] (0,0) circle [radius=0.5];
	\end{tikzpicture}
	\parbox{1.5cm}{\quad $= \delta_V$\quad}
	\begin{tikzpicture}[baseline=0]
		\draw [bgd] (0,0) circle [radius=1.0];
	\end{tikzpicture}
\end{equation*}

\begin{equation*}
	\begin{tikzpicture}[baseline=0]
		\draw [bgd] (0,0) circle [radius=1.0];
		\draw [repW] (0,0) circle [radius=0.5];
	\end{tikzpicture}
	\parbox{1.5cm}{\quad $= \delta_W$\quad}
	\begin{tikzpicture}[baseline=0]
		\draw [bgd] (0,0) circle [radius=1.0];
	\end{tikzpicture}
\end{equation*}

\subsubsection{Eigenvalues}\label{eigenvalues}
Next we introduce $\alpha$, $\beta$, $\gamma$ in the following relations.
We take $\alpha^3$ and $\gamma^2$ for later convenience.
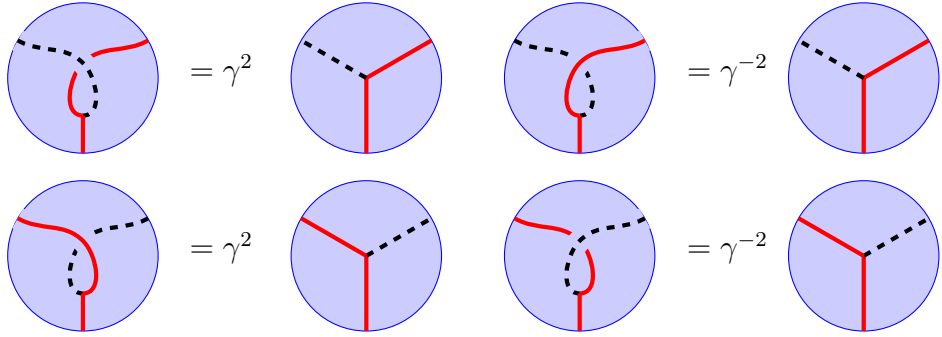
\begin{equation*}
	\begin{tikzpicture}[baseline=0]
		\draw [bgd] (0,0) circle [radius=1.0];
		\draw [repV] (1,0) to [out=180,in=45] (0,0.25) to [out=225, in=180] (0,-0.5);
		\draw [blue!20, line width=3mm] (-1,0) to [out=0,in=135] (0,0.25) to [out=-45, in=0] (0,-0.5);
		\draw [repV] (-1,0) to [out=0,in=135] (0,0.25) to [out=-45, in=0] (0,-0.5);
	\end{tikzpicture}
	\parbox{1.75cm}{\quad $=\alpha^3$\quad}
	\begin{tikzpicture}[baseline=0]
		\draw [bgd] (0,0) circle [radius=1.0];
		\draw [repV] (-1,0) -- (1,0);
	\end{tikzpicture}
	\qquad
	\begin{tikzpicture}[baseline=0]
		\draw [bgd] (0,0) circle [radius=1.0];
		\draw [repV] (-1,0) to [out=0,in=135] (0,0.25) to [out=-45, in=0] (0,-0.5);
		\draw [blue!20, line width=3mm] (1,0) to [out=180,in=45] (0,0.25) to [out=225, in=180] (0,-0.5);
		\draw [repV] (1,0) to [out=180,in=45] (0,0.25) to [out=225, in=180] (0,-0.5);
	\end{tikzpicture}
	\parbox{1.5cm}{\quad $=\alpha^{-3}$\quad}
	\begin{tikzpicture}[baseline=0]
		\draw [bgd] (0,0) circle [radius=1.0];
		\draw [repV] (-1,0) -- (1,0);
	\end{tikzpicture}
\end{equation*}

\begin{equation*}
	\begin{tikzpicture}[baseline=0]
		\draw [bgd] (0,0) circle [radius=1.0];
		\draw [repV] (30:1) to [out=210,in=45] (0,0.2) to [out=225, in=180] (0,-0.5);
		\draw [blue!20, line width=3mm] (150:1) to [out=330,in=135] (0,0.2) to [out=-45, in=0] (0,-0.5);
		\draw [repV] (150:1) to [out=330,in=135] (0,0.2) to [out=-45, in=0] (0,-0.5);
		\draw [repW] (0,-0.5) -- (0,-1);
	\end{tikzpicture}
	\parbox{1.75cm}{\quad $=\beta$\quad}
	\begin{tikzpicture}[baseline=0]
		\draw [bgd] (0,0) circle [radius=1.0];
		\draw [repW] (0,0) -- (270:1);
		\draw [repV] (0,0) -- (30:1);
		\draw [repV] (0,0) -- (150:1);
	\end{tikzpicture}
	\qquad
	\begin{tikzpicture}[baseline=0]
		\draw [bgd] (0,0) circle [radius=1.0];
		\draw [repV] (150:1) to [out=330,in=135] (0,0.2) to [out=-45, in=0] (0,-0.5);
		\draw [blue!20, line width=3mm] (30:1) to [out=210,in=45] (0,0.2) to [out=225, in=180] (0,-0.5);
		\draw [repV] (30:1) to [out=210,in=45] (0,0.2) to [out=225, in=180] (0,-0.5);
		\draw [repW] (0,-0.5) -- (0,-1);
	\end{tikzpicture}
	\parbox{1.75cm}{\quad $= \beta^{-1}$\quad}
	\begin{tikzpicture}[baseline=0]
		\draw [bgd] (0,0) circle [radius=1.0];
		\draw [repW] (0,0) -- (270:1);
		\draw [repV] (0,0) -- (30:1);
		\draw [repV] (0,0) -- (150:1);
	\end{tikzpicture}
\end{equation*}

\begin{equation*}
	\begin{tikzpicture}[baseline=0]
		\draw [bgd] (0,0) circle [radius=1.0];
		\draw [repV] (30:1) to [out=210,in=45] (0,0.2) to [out=225, in=180] (0,-0.5);
		\draw [blue!20, line width=3mm] (150:1) to [out=330,in=135] (0,0.2) to [out=-45, in=0] (0,-0.5);
		\draw [repW] (150:1) to [out=330,in=135] (0,0.2) to [out=-45, in=0] (0,-0.5);
		\draw [repV] (0,-0.5) -- (0,-1);
	\end{tikzpicture}
	\parbox{1.75cm}{\quad $=\gamma^2$\quad}
	\begin{tikzpicture}[baseline=0]
		\draw [bgd] (0,0) circle [radius=1.0];
		\draw [repV] (0,0) -- (270:1);
		\draw [repV] (0,0) -- (30:1);
		\draw [repW] (0,0) -- (150:1);
	\end{tikzpicture}
	\qquad
	\begin{tikzpicture}[baseline=0]
		\draw [bgd] (0,0) circle [radius=1.0];
		\draw [repW] (150:1) to [out=330,in=135] (0,0.2) to [out=-45, in=0] (0,-0.5);
		\draw [blue!20, line width=3mm] (30:1) to [out=210,in=45] (0,0.2) to [out=225, in=180] (0,-0.5);
		\draw [repV] (30:1) to [out=210,in=45] (0,0.2) to [out=225, in=180] (0,-0.5);
		\draw [repV] (0,-0.5) -- (0,-1);
	\end{tikzpicture}
	\parbox{1.75cm}{\quad $= \gamma^{-2}$\quad}
	\begin{tikzpicture}[baseline=0]
		\draw [bgd] (0,0) circle [radius=1.0];
		\draw [repV] (0,0) -- (270:1);
		\draw [repV] (0,0) -- (30:1);
		\draw [repW] (0,0) -- (150:1);
	\end{tikzpicture}
\end{equation*}

\begin{equation*}
	\begin{tikzpicture}[baseline=0]
		\draw [bgd] (0,0) circle [radius=1.0];
		\draw [repW] (30:1) to [out=210,in=45] (0,0.2) to [out=225, in=180] (0,-0.5);
		\draw [blue!20, line width=3mm] (150:1) to [out=330,in=135] (0,0.2) to [out=-45, in=0] (0,-0.5);
		\draw [repV] (150:1) to [out=330,in=135] (0,0.2) to [out=-45, in=0] (0,-0.5);
		\draw [repV] (0,-0.5) -- (0,-1);
	\end{tikzpicture}
	\parbox{1.75cm}{\quad $=\gamma^2$\quad}
	\begin{tikzpicture}[baseline=0]
		\draw [bgd] (0,0) circle [radius=1.0];
		\draw [repV] (0,0) -- (270:1);
		\draw [repW] (0,0) -- (30:1);
		\draw [repV] (0,0) -- (150:1);
	\end{tikzpicture}
	\qquad
	\begin{tikzpicture}[baseline=0]
		\draw [bgd] (0,0) circle [radius=1.0];
		\draw [repV] (150:1) to [out=330,in=135] (0,0.2) to [out=-45, in=0] (0,-0.5);
		\draw [blue!20, line width=3mm] (30:1) to [out=210,in=45] (0,0.2) to [out=225, in=180] (0,-0.5);
		\draw [repW] (30:1) to [out=210,in=45] (0,0.2) to [out=225, in=180] (0,-0.5);
		\draw [repV] (0,-0.5) -- (0,-1);
	\end{tikzpicture}
	\parbox{1.75cm}{\quad $=\gamma^{-2}$\quad}
	\begin{tikzpicture}[baseline=0]
		\draw [bgd] (0,0) circle [radius=1.0];
		\draw [repV] (0,0) -- (270:1);
		\draw [repW] (0,0) -- (30:1);
		\draw [repV] (0,0) -- (150:1);
	\end{tikzpicture}
\end{equation*}

\begin{lemma} The scalars $\alpha$, $\beta$, $\gamma$ satisfy $\alpha^3=\beta\gamma^2$.
\end{lemma}

\begin{proof}
We are assuming $\End(V)\cong \kk$ and $\End(W)\cong \kk$. Let $\theta_V$ and $\theta_W$ be
the values of the twists under these identifications. Then it follows from
Proposition~\ref{prop:Reigen} that
 $\alpha^3=\theta_V^2$, $\beta = \theta_V^2\theta_W^{-1}$ and $\gamma^2=\theta_W$.
 Hence 	$\alpha^3=\theta_V^2 = \beta\gamma^2$.
\end{proof}

The assumption $\End(W)\cong \kk$ implies the relation in \eqref{tadpole}.

The category $\Tcat(\alpha,\beta,\gamma,\delta_V,\delta_W,\phi,\psi)$ satisfies some additional relations.

\subsubsection{Further relations}\label{further}
In this section we make the additional assumption that $\Hom(V\otimes V,W)\cong \kk$

Define $\tau\in \kk$ by
\begin{equation}\label{triangle}
	\begin{tikzpicture}
		\draw [bgd] (0,0) circle [radius=1.0];
		\draw [repW] (90:0.5) -- (90:1);
		\draw [repV] (210:0.5) -- (210:1);
		\draw [repV] (330:0.5) -- (330:1);
		\draw [repV] (90:0.5) to [out=210, in=90] (210:0.5);
		\draw [repW] (210:0.5) to [out=330, in=210] (330:0.5);
		\draw [repV] (330:0.5) to [out=90, in=330] (90:0.5);
	\end{tikzpicture}
	\raisebox{1cm}{\quad $= \tau$\quad}
	\begin{tikzpicture}
		\draw [bgd] (0,0) circle [radius=1.0];
		\draw [repW] (0,0) -- (90:1);
		\draw [repV] (0,0) -- (210:1);
		\draw [repV] (0,0) -- (330:1);
	\end{tikzpicture}
\end{equation}

\begin{lemma} The following relations hold in $\Tcat(\alpha,\beta,\gamma,\delta_V,\delta_W,\phi,\psi)$:
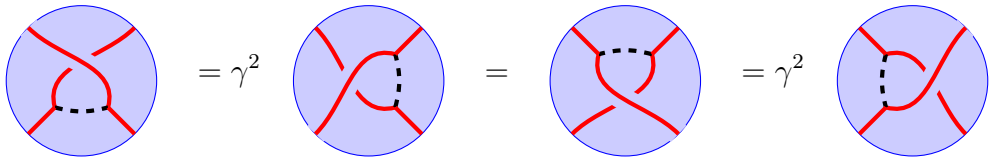
	\begin{equation*}
		\begin{tikzpicture}[baseline=0]
			\draw [bgd] (0,0) circle [radius=1.0];
			\draw [repV] (45:1) to [out=225,in=105] (225:0.5);
			\draw [blue!20, line width =3mm] (135:1) to [out=315,in=75] (315:0.5);
			\draw [repV] (135:1) to [out=315,in=75] (315:0.5);
			\draw [repW] (225:0.5) to [out=345,in=195] (315:0.5);
			\draw [repV] (225:0.5) -- (225:1);
			\draw [repV] (315:0.5) -- (315:1);
		\end{tikzpicture}
		\quad =\gamma^2\quad
		\begin{tikzpicture}[baseline=0,rotate=90]
			\draw [bgd] (0,0) circle [radius=1.0];
			\draw [repV] (45:1) to [out=225,in=105] (225:0.5);
			\draw [blue!20, line width =3mm] (135:1) to [out=315,in=75] (315:0.5);
			\draw [repV] (135:1) to [out=315,in=75] (315:0.5);
			\draw [repW] (225:0.5) to [out=345,in=195] (315:0.5);
			\draw [repV] (225:0.5) -- (225:1);
			\draw [repV] (315:0.5) -- (315:1);
		\end{tikzpicture}
		\quad =\quad
		\begin{tikzpicture}[baseline=0,rotate=180]
			\draw [bgd] (0,0) circle [radius=1.0];
			\draw [repV] (45:1) to [out=225,in=105] (225:0.5);
			\draw [blue!20, line width =3mm] (135:1) to [out=315,in=75] (315:0.5);
			\draw [repV] (135:1) to [out=315,in=75] (315:0.5);
			\draw [repW] (225:0.5) to [out=345,in=195] (315:0.5);
			\draw [repV] (225:0.5) -- (225:1);
			\draw [repV] (315:0.5) -- (315:1);
		\end{tikzpicture}
		\quad =\gamma^2\quad
		\begin{tikzpicture}[baseline=0,rotate=270]
			\draw [bgd] (0,0) circle [radius=1.0];
			\draw [repV] (45:1) to [out=225,in=105] (225:0.5);
			\draw [blue!20, line width =3mm] (135:1) to [out=315,in=75] (315:0.5);
			\draw [repV] (135:1) to [out=315,in=75] (315:0.5);
			\draw [repW] (225:0.5) to [out=345,in=195] (315:0.5);
			\draw [repV] (225:0.5) -- (225:1);
			\draw [repV] (315:0.5) -- (315:1);
		\end{tikzpicture}
	\end{equation*}
\end{lemma}

\begin{defn} Define the diagram
	
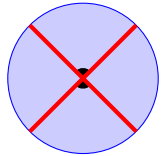
\begin{equation*}
		\begin{tikzpicture}[baseline=0]
			\draw [bgd] (0,0) circle [radius=1.0];
			\posdot{0,0};
			\draw[repV] (45:1) -- (225:1);
			\draw[repV]	(135:1) -- (315:1);
		\end{tikzpicture}
	\end{equation*}
to be any of the following

\begin{equation*}
	\gamma^{-1}\:
	\begin{tikzpicture}[baseline=0]
		\draw [bgd] (0,0) circle [radius=1.0];
		\draw [repV] (45:1) to [out=225,in=105] (225:0.5);
		\draw [blue!20, line width =3mm] (135:1) to [out=315,in=75] (315:0.5);
		\draw [repV] (135:1) to [out=315,in=75] (315:0.5);
		\draw [repW] (225:0.5) to [out=345,in=195] (315:0.5);
		\draw [repV] (225:0.5) -- (225:1);
		\draw [repV] (315:0.5) -- (315:1);
	\end{tikzpicture}
\: =\gamma\:
	\begin{tikzpicture}[baseline=0,rotate=90]
		\draw [bgd] (0,0) circle [radius=1.0];
		\draw [repV] (45:1) to [out=225,in=105] (225:0.5);
		\draw [blue!20, line width =3mm] (135:1) to [out=315,in=75] (315:0.5);
		\draw [repV] (135:1) to [out=315,in=75] (315:0.5);
		\draw [repW] (225:0.5) to [out=345,in=195] (315:0.5);
		\draw [repV] (225:0.5) -- (225:1);
		\draw [repV] (315:0.5) -- (315:1);
	\end{tikzpicture}
	 =\gamma^{-1}\:
	\begin{tikzpicture}[baseline=0,rotate=180]
		\draw [bgd] (0,0) circle [radius=1.0];
		\draw [repV] (45:1) to [out=225,in=105] (225:0.5);
		\draw [blue!20, line width =3mm] (135:1) to [out=315,in=75] (315:0.5);
		\draw [repV] (135:1) to [out=315,in=75] (315:0.5);
		\draw [repW] (225:0.5) to [out=345,in=195] (315:0.5);
		\draw [repV] (225:0.5) -- (225:1);
		\draw [repV] (315:0.5) -- (315:1);
	\end{tikzpicture}
	\: =\gamma\:
	\begin{tikzpicture}[baseline=0,rotate=270]
		\draw [bgd] (0,0) circle [radius=1.0];
		\draw [repV] (45:1) to [out=225,in=105] (225:0.5);
		\draw [blue!20, line width =3mm] (135:1) to [out=315,in=75] (315:0.5);
		\draw [repV] (135:1) to [out=315,in=75] (315:0.5);
		\draw [repW] (225:0.5) to [out=345,in=195] (315:0.5);
		\draw [repV] (225:0.5) -- (225:1);
		\draw [repV] (315:0.5) -- (315:1);
	\end{tikzpicture}
\end{equation*}
\end{defn}

\subsubsection{Bigons}\label{bubbles}
Define $\phi,\psi\in\kk^*$ by
\begin{equation}\label{eq:bubbleV}
\begin{tikzpicture}[baseline=0]
	\draw [bgd] (0,0) circle [radius=1.0];
	\draw [repW] (0,0.5) to [out=210, in=90] (-0.4,0) to [out=270, in=150] (0,-0.5);
	\draw [repV] (0,0.5) to [out=330, in=90] (0.4,0) to [out=270, in=30](0,-0.5);
	\draw [repV] (0,0.5) -- (0,1);
	\draw [repV] (0,-0.5) -- (0,-1);
\end{tikzpicture}
\parbox[b]{1.5cm}{\quad$=\phi$}
\begin{tikzpicture}[baseline=0]
	\draw [bgd] (0,0) circle [radius=1.0];
	\draw [repV] (0,1) -- (0,-1);
\end{tikzpicture}
\end{equation}

\begin{equation}\label{eq:bubbleW}
\begin{tikzpicture}[baseline=0]
	\draw [bgd] (0,0) circle [radius=1.0];
	\draw [repV] (0,0.5) to [out=210, in=90] (-0.4,0) to [out=270, in=150] (0,-0.5);
	\draw [repV] (0,0.5) to [out=330, in=90] (0.4,0) to [out=270, in=30](0,-0.5);
	\draw [repW] (0,0.5) -- (0,1);
	\draw [repW] (0,-0.5) -- (0,-1);
\end{tikzpicture}
\parbox[b]{1.5cm}{\quad$= \psi$}
\begin{tikzpicture}[baseline=0]
	\draw [bgd] (0,0) circle [radius=1.0];
	\draw [repW] (0,1) -- (0,-1);
\end{tikzpicture}
\end{equation}

\begin{lemma}\label{theta}
These satisfy $\phi\delta_V = \psi\delta_W$.
\end{lemma}
\begin{proof}	
\begin{equation}
	\parbox[b]{1.5cm}{$\phi\delta_V=$}
	\begin{tikzpicture}[baseline=0]
		\draw [bgd] (0.5,0) circle [radius=1.5];
		\draw [repW] (0,0.5) to [out=210, in=90] (-0.4,0) to [out=270, in=150] (0,-0.5);
		\draw [repV] (0,0.5) to [out=330, in=90] (0.4,0) to [out=270, in=30](0,-0.5);
		\draw [repV] (0,0.5) arc(180:0:0.5) -- (1,-0.5) arc(360:180:0.5);
	\end{tikzpicture}
	\parbox[b]{1.5cm}{\quad$=$}
	\begin{tikzpicture}[baseline=0]
		\draw [bgd] (0.5,0) circle [radius=1.5];
		\draw [repV] (0,0.5) to [out=210, in=90] (-0.4,0) to [out=270, in=150] (0,-0.5);
		\draw [repV] (0,0.5) to [out=330, in=90] (0.4,0) to [out=270, in=30](0,-0.5);
		\draw [repW] (0,0.5) arc(180:0:0.5) -- (1,-0.5) arc(360:180:0.5);
	\end{tikzpicture}
	\parbox[b]{2cm}{\quad$=\psi\delta_W$}
\end{equation}
\end{proof}

\section{Endomorphism algebra}\label{sec:relations}
In this section we study the endomorphism algebra
\begin{equation*}
	A(2) = \End\left(
\begin{tikzpicture}[baseline=0]
	\draw[bgd] (0,-0.4) rectangle (1,0.6);
	\draw[repV] (0.25,0.6) -- (0.25,-0.4);
	\draw[repV] (0.75,0.6) -- (0.75,-0.4);
\end{tikzpicture}
	\right) = \End(V\otimes V)
\end{equation*}

This has additional structure, the rotation map, which is a linear map
$\Rot\colon A(2)\to A(2)$ and the trace map, which is a linear map
$\Tr\colon A(2)\to \kk$.

\begin{defn}
The rotation map is the linear map, $\Rot$, given by
\begin{equation}\label{fig:rot}
	\Rot\colon
	\begin{tikzpicture}[baseline=0]
		\draw[bgd] (-0.1,-0.8) rectangle (1.1,1);
		\draw[repV] (0.25,1) -- (0.25,-0.8);
		\draw[repV] (0.75,1) -- (0.75,-0.8);
		\filldraw [fill=Peach, draw=black] (0.1,-0.2) rectangle (0.9,0.4);
	\end{tikzpicture}
\:\mapsto\:
\begin{tikzpicture}[baseline=0]
	\draw[bgd] (-0.5,-0.8) rectangle (1.5,1);
	\filldraw [fill=Peach, draw=black] (0.1,-0.2) rectangle (0.9,0.4);
	\draw[repV] (0.25,0.4) arc (0:180:0.25) -- (-0.25,-0.8);
	\draw[repV] (0.75,-0.2) arc (180:360:0.25) -- (1.25,1);
	\draw[repV] (0.75,0.4) -- (0.75,1);
	\draw[repV] (0.25,-0.2) -- (0.25,-0.8);
\end{tikzpicture}
\end{equation}
\end{defn}

\begin{defn}\label{defn:tr}
	The trace map is the linear map, $\Tr$, given by
	\begin{equation}\label{fig:tr}
		\Tr\colon
		\begin{tikzpicture}[baseline=0]
			\draw[bgd] (-0.1,-0.8) rectangle (1.1,1);
			\draw[repV] (0.25,1) -- (0.25,-0.8);
			\draw[repV] (0.75,1) -- (0.75,-0.8);
			\filldraw [fill=Peach, draw=black] (0.1,-0.2) rectangle (0.9,0.4);
		\end{tikzpicture}
		\:\mapsto\:
		\begin{tikzpicture}[baseline=0]
			\draw[bgd] (-0.1,-0.8) rectangle (1.5,1);
			\draw[repV] (0.25,1) -- (0.25,-0.8);
			\filldraw [fill=Peach, draw=black] (0.1,-0.2) rectangle (0.9,0.4);
			\draw[repV] (0.75,-0.2) arc (180:360:0.25) -- (1.25,0.4) arc (0:180:0.25);
		\end{tikzpicture}
	\end{equation}
\end{defn}

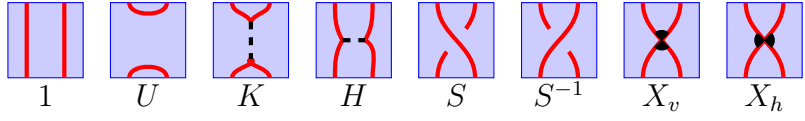
\begin{figure}[ht]
\begin{tabular}{cccccccc}
\begin{tikzpicture}[baseline=0] 
	\draw[bgd] (0,-0.4) rectangle (1,0.6);
	\draw[repV] (0.25,0.6) -- (0.25,-0.4);
	\draw[repV] (0.75,0.6) -- (0.75,-0.4);
\end{tikzpicture}
&
\begin{tikzpicture}[baseline=0] 
	\draw[bgd] (0,-0.4) rectangle (1,0.6);
	\draw[repV] (0.25,0.6) to[out=270,in=270] (0.75,0.6);
	\draw[repV] (0.25,-0.4) to[out=90,in=90] (0.75,-0.4);
\end{tikzpicture}
&\begin{tikzpicture}[baseline=0] 
	\draw[bgd] (0,-0.4) rectangle (1,0.6);
	\draw[repV] (0.25,0.6) to[out=270,in=120] (0.5,0.35);
	\draw[repV] (0.75,0.6) to[out=270,in=60] (0.5,0.35);
	\draw[repV] (0.25,-0.4) to[out=90,in=330] (0.5,-0.15);
	\draw[repV] (0.75,-0.4) to[out=90,in=210] (0.5,-0.15);
	\draw[repW] (0.5,0.35) -- (0.5,-0.15);
\end{tikzpicture}
&\begin{tikzpicture}[baseline=0] 
	\draw[bgd] (0,-0.4) rectangle (1,0.6);
	\draw[repV] (0.25,0.6) to[out=270,in=120] (0.35,0.1);
	\draw[repV] (0.75,0.6) to[out=270,in=60] (0.65,0.1);
	\draw[repV] (0.25,-0.4) to[out=90,in=240] (0.35,0.1);
	\draw[repV] (0.75,-0.4) to[out=90,in=330] (0.65,0.1);
	\draw[repW] (0.35,0.1) -- (0.65,0.1);
\end{tikzpicture}
&\begin{tikzpicture}[baseline=0]
	\draw[bgd] (0,-0.4) rectangle (1,0.6);
	\draw[repV] (0.75,0.6) to[out=270,in=90] (0.25,-0.4);
	\fill[blue!20] (0.5,0.1) circle (0.2);
	\draw[repV] (0.25,0.6) to[out=270,in=90] (0.75,-0.4);
\end{tikzpicture}
&\begin{tikzpicture}[baseline=0]
	\draw[bgd] (0,-0.4) rectangle (1,0.6);
	\draw[repV] (0.25,0.6) to[out=270,in=90] (0.75,-0.4);
	\fill[blue!20] (0.5,0.1) circle (0.2);
	\draw[repV] (0.75,0.6) to[out=270,in=90] (0.25,-0.4);
\end{tikzpicture}
&\begin{tikzpicture}[baseline=0]
	\draw[bgd] (0,-0.4) rectangle (1,0.6);
	\draw[repV] (0.25,0.6) to[out=270,in=90] (0.75,-0.4);
	\draw[repV] (0.75,0.6) to[out=270,in=90] (0.25,-0.4);
	\posdot{0.5,0.1}
\end{tikzpicture}
&\begin{tikzpicture}[baseline=0]
	\draw[bgd] (0,-0.4) rectangle (1,0.6);
	\draw[repV] (0.25,0.6) to[out=270,in=90] (0.75,-0.4);
	\draw[repV] (0.75,0.6) to[out=270,in=90] (0.25,-0.4);
	\negdot{0.5,0.1}
\end{tikzpicture} \\
1 & $U$ & $K$ & $H$ & $S$ & $S^{-1}$ & $X_v$ & $X_h$
\end{tabular}
\caption{Elements of two point algebra}\label{fig:elements}
\end{figure}

We name the elements as in Figure~\ref{fig:elements}.

\begin{assume}\label{assumption}
The set $\{1,U,K,H,S+S^{-1}\}$ is a basis of $A(2)$.
\end{assume}
This implies that rotation has order two and acts by
$1\leftrightarrow U$, $K\leftrightarrow H$,
$S\leftrightarrow S^{-1}$, $X_h\leftrightarrow X_v$.

Furthermore, the $+1$-eigenspace has basis $\{(1+U),(H+K),(S+S^{-1})\}$
and the $-1$-eigenspace has basis by $\{(1-U),(H-K)\}$.

This implies that we have a relation of the form
\begin{equation}\label{eq:skein}
	S - S^{-1} = z(1-U)+w(H-K)
\end{equation}
This relation is called the \Dfn{skein relation}.

It follows that the set
\begin{equation}\label{basis}
	\bB=\{1,U,K,H,S\}
\end{equation}
is also a basis.

The trace map is given on the basis \eqref{basis} by
\begin{equation*}
	\Tr(1)=\delta_V \quad
	\Tr(U)=1 \quad
	\Tr(K)=\phi \quad
	\Tr(H)=0 \quad
	\Tr(S)=\alpha^{-3}
\end{equation*}

The rotation map is a linear map $\Rot\colon A(2)\to A(2)$. This satisfies $\Rot^2=1$.

The matrix of $\Rot$ with respect to the basis $\bB$ is
\begin{equation}\label{rotmatrix}
	\Rot = \begin{pmatrix}
		0 & 1 & 0 & 0 & 0 \\
		1 & 0 & 0 & 0 & 0 \\
		0 & 0 & 0 & 1 & 0 \\
		0 & 0 & 1 & 0 & 0 \\
		-z & z & -w & w & 1 \\
	\end{pmatrix}
\end{equation}

Using $\Rot$ and the maps $S^{\pm 1}$ and the trace map means that we have an element of
$D_3(t,u,v)$ for each rational link. Later we can explain that these values interpolate between the link polynomials
in the sense that these values specialise to the link polynomials.

\begin{lemma} These satisfy $(S\;\Rot)^3= \alpha^3 = (\Rot\; S)^3 $.
\end{lemma}
\begin{proof}
The proof is based on \cite{Conway1970}. The difference is that Conway uses unoriented links
and tangles and has the Reidemeister I move whereas we use unoriented framed links and tangles.
\end{proof}

\begin{lemma}
Put $\sigma_1=S$ and $\sigma_2=\Rot\;S\;\Rot$. Then it follows that these
satisfy the braid relation.
\begin{equation*}
	\sigma_1\sigma_2\sigma_1 = \alpha^3\Rot = \sigma_2\sigma_1\sigma_2
\end{equation*}
\end{lemma}

\begin{proof}
	\begin{equation*}
\sigma_1\sigma_2\sigma_1 = (S\;\Rot)^3\Rot =  \alpha^3\Rot = (\Rot\; S)^3\Rot	
	\end{equation*}
\end{proof}
\subsection{Structure constants}
\begin{thm}\label{thm:constants}
The structure constants of $A(2)$ with respect to
the basis $\mathcal{B}$ are rational functions
of $\alpha,\beta,\gamma,\delta_V,\delta_W,\phi,\psi$.
\end{thm}

This is proved in the rest of this section by deriving
explicit rational functions for the structure constants.

The initial relations in \S~\ref{initial} give
\begin{equation*}
U^2=\delta_V U \qquad U\,K = 0 = K\,U	
\end{equation*}

The eigenvalue relations in \S~\ref{eigenvalues} give
\begin{gather*}
	U\,S^{\pm 1}=\alpha^{\pm 3} U = S^{\pm 1}\,U \\
	K\,S^{\pm 1}=\beta^{\pm 1} K = S^{\pm 1}\,K	
\end{gather*}

The relations in \S~\ref{further} give
\begin{gather*}
	S\,H = \gamma X_h = H\,S \\
	S^{-1}\,H = \gamma^{-1}X_v = H\,S^{-1}
\end{gather*}

The bigon relations in \S~\ref{bubbles} give
\begin{equation*}
	U\,H=\phi\, U = H\,U \qquad K^2 = \psi\, K
\end{equation*}

The multiplication table of $A(2)$ with respect to the basis $\bB$ is shown in Figure~\ref{fig:table}. In order
to complete this table it remains to determine the
coefficient $\tau$ (defined in \eqref{triangle}) and to express $X_v,S^2,H^2$ in terms of the spanning set, $\mathcal{B}$, defined.in \eqref{basis}. These expressions are
given in Proposition~\ref{HSrelation}, Proposition~\ref{digonrelation}
and Proposition~\ref{squarerelation}.

\begin{figure}[ht]
\begin{tabular}{c|cccc}
& $U$ & $K$ & $H$ & $S$ \\ \hline
$U$ & $\delta_V U$ & 0 & $\phi U$ & $\alpha^3 U$ \\
$K$ & 0 & $\psi K$ & $\tau K$ & $\beta K$ \\
$H$ & $\phi U$ & $\tau K$ & $H^2$ & $\gamma X_h$ \\
$S$ & $\alpha^3 U$ & $\beta K$ & $\gamma X_h$ & $S^2$ \\
\end{tabular}
\caption{Multiplication table}\label{fig:table}
\end{figure}

\begin{prop}\label{HSrelation} The following relation holds:
	\begin{equation*}
X_h = a_I(1+\alpha^{2} U+\alpha S) -\alpha\gamma^{-1} H
-\alpha^{-1}\gamma K
	\end{equation*}
\end{prop}
\begin{proof}
We have a relation with unknown coefficients
\begin{equation}\label{eq:Xh}
	X_h = a_I + a_U U + a_K K + a_H H + a_S S
\end{equation}
Multiply by $S^{-1}$ to get
\begin{equation*}
	\gamma^{-1} H = a_I S^{-1}  + \alpha^{-3} a_U U + \beta^{-1} a_K K + \gamma^{-1} a_H X_v + a_S
\end{equation*}
Rotating gives
\begin{equation*}
	\gamma^{-1} K = a_I S  + \alpha^{-3} a_U + \beta^{-1} a_K H + \gamma^{-1} a_H X_h + a_S U
\end{equation*}
Solve for $a_HX_h$ to get
\begin{equation}\label{eq:HS}
	a_H X_h = -\alpha^{-1}\gamma a_U -\gamma a_S U + K -\beta^{-1}\gamma a_K H - \gamma a_I S
\end{equation}
Multiply \eqref{eq:Xh} by $a_H$ to get a second equation for $a_H X_h$.
By the assumption that the set $\bB$ is a basis the individual coefficients
can be compared.

Comparing coefficients of $K$ gives $a_K a_H = 1$.
Comparing coefficients of $H$ gives $a_H^2 = -\beta^{-1}\gamma a_K$.

Solving this pair of equations gives
\begin{equation*}
a_K = -(\alpha\gamma^{-1})\qquad a_H = -(\alpha^{-1}\gamma)
\end{equation*}

The remaining three equations give
\begin{equation*}
	a_I = \alpha^{-2}a_U,\qquad
	a_U = \alpha a_S,\qquad
	a_S = \alpha a_I
\end{equation*}
\end{proof}

\begin{cor}\label{HSinv}
	\begin{equation*}
	X_v = a_I(\alpha^{2}+U+\alpha S^{-1}) -\alpha\gamma^{-1} H
	-\alpha^{-1}\gamma K 
\end{equation*}
\end{cor}
\begin{proof}
	This is given by rotating the expression for
	$X_h$ in Proposition~\ref{HSrelation}.
\end{proof}

\begin{lemma} If $1+\alpha^2\delta_V+\alpha^4\ne 0$
	\begin{equation*}
a_I = \alpha^{-1}\phi\left(\frac{\alpha^2\gamma^{-1}+\alpha^{-2}\gamma}%
{\alpha^{-2}+\delta_V+\alpha^2}\right)
\end{equation*}
	\end{lemma}
\begin{proof}
Multiplying Proposition~\ref{HSrelation} by $U$ gives
\begin{equation*}
	\alpha^3\gamma^{-1}\phi = a_I(1+\alpha^2\delta_V+\alpha^4) -  \alpha^{-1}\gamma\phi
\end{equation*}
Now solve for $a_I$.
\end{proof}

\begin{lemma} If $\beta\gamma^{-1}+\alpha^{-1}\gamma\ne 0$
	\begin{equation*}
	\tau = \frac{\alpha a_I(\alpha^{-2}\gamma+\alpha^2\gamma^{-1})-\psi}%
{(\alpha^{-2}\gamma^2+\alpha^2\gamma^{-2})}
	\end{equation*}
\end{lemma}
\begin{proof}	
	Multiplying Proposition~\ref{HSrelation} by $K$ gives
\begin{equation*}
	\beta\gamma^{-1}\tau = a_I(1+\alpha\beta)-\alpha\gamma^{-1}\psi -\alpha^{-1}\gamma\tau
\end{equation*}
	Solve for $\tau$ to get
	\begin{equation*}
	\tau = \frac{a_I(1+\alpha\beta)-\alpha\gamma^{-1}\psi}{\beta\gamma^{-1}+\alpha^{-1}\gamma}
\end{equation*}
	Substitute $\beta=\alpha^3\gamma^{-2}$ to get
\begin{equation*}
	\tau = \frac{\alpha a_I(\alpha^{-2}\gamma+\alpha^2\gamma^{-1})-\psi}%
	{(\alpha^{-2}\gamma^2+\alpha^2\gamma^{-2})}
\end{equation*}	
\end{proof}

\begin{prop}\label{skeinrelation} Assume $(\tau-\psi)\delta_V+\phi+\psi-\tau\ne 0$.
	Then the coefficients in the skein relation, \eqref{eq:skein}, are given by
\begin{align*}
	z &= \frac{(\beta-\beta^{-1})\phi+(\alpha^3-\alpha^{-3})(\psi-\tau)}{(\tau-\psi)\delta_V+\phi+\psi-\tau} \\
	w &= \frac{(\beta-\beta^{-1})(\delta_V-1)+(\alpha^3-\alpha^{-3})}{(\tau-\psi)\delta_V+\phi+\psi-\tau}
\end{align*}
\end{prop}
\begin{proof}
Multiplying the skein relation by $U$ and by $K$ gives
\begin{equation*}
	\alpha^3-\alpha^{-3} = z(1-\delta_V) + w \phi
	\qquad
	\beta-\beta^{-1} = z + w(\tau-\psi)
\end{equation*}
Solving this pair of equations for $z$ and $w$ gives the result.
\end{proof}

The next Lemma is a variation of the skein relation
\begin{lemma}\label{lem:skX}
\begin{equation*}
	X_h - X_v = (\alpha^{-1}-\alpha+z)(1-U)+(\alpha^{-1}\gamma-\alpha\gamma^{-1}+w)(H-K)
\end{equation*}	
\end{lemma}
\begin{proof} Use the equations for $X$ to get
	\begin{equation*}
X_h - X_v = (\alpha^{-1}-\alpha)(1-U)+(\alpha^{-1}\gamma-\alpha\gamma^{-1})(H-K) + S - S^{-1}		
	\end{equation*}
Now use the skein relation.
\end{proof}
\begin{cor}
	\begin{gather*}
		(\alpha^3\gamma^{-1}-\alpha^{-3}\gamma)\phi = (\alpha^{-1}-\alpha+z)(1-\delta_V)+(\alpha^{-1}\gamma-\alpha\gamma^{-1}+w)\phi\\
		(\beta\gamma^{-1}-\beta^{-1}\gamma)\tau = (\alpha^{-1}-\alpha+z)+(\alpha\gamma^{-1}-\alpha^{-1}\gamma+w)(\tau-\psi)
	\end{gather*}
\end{cor}
\begin{proof}
Multiply Lemma~\ref{lem:skX} by $U$ and by $K$ to get the equations.
\end{proof}

\begin{prop}\label{digonrelation}
\begin{equation*}
	S^2 = 1 + \alpha^{-1}\gamma w + \alpha(\gamma w - \alpha^2 z) U - (\alpha+\beta) w K
	- \alpha^{-1}\gamma^2 w H + (z+\gamma w) S
\end{equation*}
\end{prop}
\begin{proof}
Multiply the skein relation, \eqref{eq:skein}, by $S$ and solve for $S^2$. This gives
	\begin{equation*}
	S^2 = 1 + z S - z\alpha^3 U + w\gamma X_h -w\beta K
\end{equation*}
Now substitute for $X_h$ using Proposition~\ref{HSrelation}.
\end{proof}
\begin{prop}\label{squarerelation}
	\begin{multline*}
		wH^2 = (\alpha^{-1}\gamma-\alpha\gamma^{-1}-\gamma^{-1}z) 
		+(\alpha\gamma-\alpha^{-1}\gamma^{-1}+\gamma^{-1}z+z\phi)\;U \\
		-(\alpha^{-1}\gamma^2-\alpha\gamma^{-2}+\gamma^{-1}w+z)\;H 
		+(\alpha^{-1}-\alpha+\gamma^{-1}w+w\tau)\;K 
		+ (\gamma+\gamma^{-1})\;S
	\end{multline*}
\end{prop}

\begin{proof}
	Multiply the skein relation, \eqref{eq:skein}, by $H$ to get
\begin{align*}
	H\,S - H\,S^{-1}  &= zH -z\psi U + w H^2 - w\tau K \\
		\gamma X_h - \gamma^{-1}X_v  &= zH -z\psi U + w H^2 - w\tau K
	\end{align*}
Substitute for $X_h$ using Proposition~\ref{HSrelation} and for $X_v$ using Corollary~\ref{HSinv}.
\begin{multline*}
	\gamma X_h - \gamma^{-1}X_v  = (\alpha^{-1}\gamma-\alpha\gamma^{-1})
	+ (\alpha\gamma-\alpha^{-1}\gamma^{-1})U \\
	- (\alpha^{-1}\gamma^2-\alpha\gamma^{-2})H
	- (\alpha-\alpha^{-1}) K + \gamma S - \gamma^{-1} S^{-1}
\end{multline*}
\end{proof}

\section{The diagram category}\label{sec:quantum}
The previous section gives a uniform construction of relations and evaluation functors
for each example. In this section we give a universal construction. The examples are
then given by specialising.

It would be cumbersome to work with the variables and equations in \S~\ref{sec:Tdef}.
In this section we define a commutative ring of scalars and give elements which
satisfy the equations. These relations are given by a direct calculation and we
refer the interested reader to the notebook for the verifications of these equations.

The specialisations needed for the evaluation functors take a simple form and this
was the original motivation for the construction of this ring.

\subsection{Scalars}\label{scalars}
Let $J$ denote the maximal ideal of $\bQ[q,q^{-1}]$ generated by $q-1$.  Equivalently, $J$ is the kernel of the $\bQ$-algebra homomorphism $\bQ[q,q^{-1}] \to \bQ$, $q \mapsto 1$.  Let
\begin{equation} \label{Adef}
	A := \left\{ \frac{a}{b} : a,b \in \bQ[q,q^{-1}],\ b \notin J \right\} \subseteq \bQ(q)
\end{equation}
be the localization of $\bQ[q,q^{-1}]$ at the ideal $J$.  Since $J$ is also the maximal ideal of $\bQ[q,q^{-1}]$ generated by $q^{-1}-1$, the map
\[
\bar{\ } \colon \bQ[q,q^{-1}] \to \bQ[q,q^{-1}],\quad \overline{q^{\pm 1}} = q^{\mp 1},
\]
induces a $\bQ$-linear involution of $A$.  The group of units of $A$ is
\[
A^\times = \left\{ \frac{a}{b} : a,b \in \bQ[q,q^{-1}],\ a,b \notin J \right\},
\]
and $A$ is a local ring with unique maximal ideal $(q-1)A$. In particular, $A$ is a principal ideal domain.

We write $\bf{t}$ for a sequence of indeterminates, ${\bf t}=(t_1,t_2,\dotsc,t_r)$,
and $\bf{I}$ for a sequence of integers, ${\bf I}=(i_1,i_2,\dotsc,i_r)$. We write
$\bf{t}^{\bf{I}}$ for the monomial $t_1^{i_1}t_2^{i_2}\dotsb t_r^{i_r}$.

Denote the Laurent polynomial ring by $\bQ[\bf{t}^{\pm 1}]$ and the field of rational functions by $\bQ(\bf{t})$. The bar involution is defined on both rings by
$t_i^{\pm 1}\mapsto t_i^{\mp 1}$ so the inclusion $\bQ[\bf{t}^{\pm 1}]\to \bQ(\bf{t})$
is a homomorphism of rings with involution.

For a sequence $\bf{I}$, put  $[\bf{I}]=\bf{t}^{\bf{I}}-\bf{t}^{-\bf{I}}\in \bQ[\bf{t}^{\pm 1}]$. Invert all elements $[\bf{I}]$
to get the localisation $\{\bf{I}\}^{-1}\bQ[\bf{t}^{\pm 1}]$. Define $D(\bf{t})$ to be the
$\bQ[\bf{t}^{\pm 1}]$-subalgebra of $\{\bf{I}\}^{-1}\bQ[\bf{t}^{\pm 1}]$ generated by elements $\frac{[\bf{I}]}{[\bf{J}]}$. The bar involution extends to this ring by
$[\bf{I}]\mapsto -[\bf{I}]$.

Note that
\begin{equation*}
	\frac{[2\bf{I}]}{[\bf{I}]} = \bf{t}^{\bf{I}}+\bf{t}^{-\bf{I}}
\end{equation*}

Let $\bQ^{(0)}(\bf{m})$ be the field of homogeneous rational functions of degree 0.
This is a ring with involution the identity homomorphism.
\begin{defn}
The homomorphism of rings with involution $\ev\colon D(\bf{t}) \to \bQ^{(0)}(\bf{m})$ is defined by
\begin{equation*}
	t_i \mapsto 1 \qquad \frac{[I]}{[J]} \mapsto \frac{i_1m_1+\dotsb +i_km_k}{j_1m_1+\dotsb +j_km_k}
\end{equation*}
\end{defn}

A matrix $P$ gives a homomorphism $\phi(P)\colon \bQ(\bf{t}) \to \bQ(\bf{t}')$.
This homomorphism is given by $t_i\mapsto \bf{t}^{P_i}$ where $P_i$ is a row/column of $P$.
This restricts to a homomorphism $\phi(P)\colon D(\bf{t})\to D(\bf{t}')$. Similarly,
there is a homomorphism $\phi(P)\colon \bQ^{(0)}(\bf{m}) \to \bQ^{(0)}(\bf{m}')$.

The following diagrams commute:
\begin{equation*}
	\begin{CD}
		D(\bf{t}) @>{\phi(P)}>> D(\bf{t}') \\
		@V{\ev}VV @VV{\ev}V \\
		\bQ^{(0)}(\bf{m}) @>>{\phi(P)}> \bQ^{(0)}(\bf{m}')
	\end{CD}\qquad
	\begin{CD}
		D(\bf{t}) @>{\phi(P)}>> D(\bf{t}) \\
		@VVV @VVV \\
		\bQ(\bf{t}) @>>{\phi(P)}>\bQ(\bf{t})
	\end{CD}
\end{equation*}

Note that composition of the homomorphisms $\phi(M)$ corresponds to matrix multiplication.
More precisely $\phi(M) \circ \phi(M') = \phi(MM')$.

The simplest case is when $\bf{t}$ is a sequence of length one. Take the single entry
to be $q$, Then the ring $\bQ[\bf{t}^{\pm 1}]$ is $\bQ[q,q^{-1}]$ and the field  $\bQ(\bf{t})$ is $\bQ(q)$. For $\bf{I}$ and $\bf{J}$ sequences of length one with entries
$i$ and $j$  we have $\frac{[\bf{I}]}{[\bf{J}]} = \frac{[i]}{[j]}$ where $[i]$
and $[j]$ are the usual quantum integers. There is an inclusion $D(q)\to A$
compatible with the bar involution. For $\bf{m}$ a sequence of length one we identify
$\bQ^{(0)}$ with $\bQ$. Then $\ev\colon D(q)\to \bQ$ is just evaluation at $q=1$.

Let $P$ be a matrix with one row/column. Then we have the homomorphism
$\phi(P)\colon \bQ(\bf{t})\to \bQ(q)$ given by $t_i\mapsto q^{p_i}$
and the homomorphism $\phi(P)\colon \bQ^{(0)}(\bf{m}) \to \bQ$ given by
$m_i\mapsto p_i$. Then the following diagram commutes
\begin{equation*}
	\begin{CD}
D(\bf{t}) @>{\phi(P)}>> A \\
@V{\ev}VV @VV{q\to 1}V \\
\bQ^{(0)}(\bf{m}) @>>{\phi(P)}> \bQ	
	\end{CD}
\end{equation*}

In this paper we will use the cases ${\bf t}=(t,u,v)$, ${\bf t}=(q,z)$ and ${\bf t}=(q)$.

We will use a variation of the notation $[{\bf I}]$ and use 
\begin{equation*}
	[rk+s] = \frac{q^sz^r-q^{-s}z^{-r}}{q-q^{-1}}v 
\end{equation*}
In particular, $[s]$ is the usual quantum integer.

\subsection{Relations}\label{algebra}
The values of $\alpha$,$\beta$,$\gamma$ are
\begin{equation}\label{eq:twists}
	\alpha = t^2 u^2  v^2\qquad
	\beta = v^{2}\qquad
	\gamma = t^3 u^3 v^2
\end{equation}

It is clear that these satisfy $\alpha^3=\beta\gamma^2$.

The quantum dimension of $V$ is
\begin{equation}\label{eq:deltaV} \delta_V =
(-1)\frac{[2,2,2][1,2,1][4,2,4]}%
{[1,1,1][2,4,2][2,1,2]}
\frac{[2,3,3][3,3,2]}{[1,0,0][0,0,1]}
\end{equation}

The quantum dimension of $W$ is
\begin{equation}\label{eq:deltaW} \delta_W =
\frac{[3,3,0][2,2,1][1,2,0][2,2,2][3,3,3][2,3,2]
	[2,3,3][1,2,1][4,2,4][4,2,2]}%
{[0,0,2][0,0,1][1,0,-1][1,0,0][1,1,0][2,4,0]
	[1,1,1][2,4,2][2,1,2][2,1,1]}
\end{equation}

Our choice of $\psi$ and $\phi$ is
\begin{align}
	\psi &= \frac{[1,0,-1][0,0,2][1,1,0][3,3,2][2,4,0][1,2,1]}%
	{[1,0,0][0,0,1][3,3,0][2,2,1][1,2,0][2,4,2]} \label{eq:psi} \\
\phi &= (-1)\frac{[3,3,3][2,3,2][1,2,1][4,2,2]}%
{[0,0,1][1,0,0][2,4,2][2,1,1]} \label{eq:phi}
\end{align}
By a direct calculation, these satisfy $\phi\;\delta_V = \psi\;\delta_W$.

The category $\Dcat$ is a $D(t,u,v)$-linear category. Take $\Tcat$ to be the 
category constructed above by taking $\kk=D(t,u,v)$ and the parameters as above.
Let $\mathcal{I}$ be the tensor ideal generated by the relations. Then
$\Dcat = \Tcat/I$ is a $D(t,u,v)$-linear ribbon category.

The following sections give the relations where the coefficients are now taken
to be elements of $D(t,u,v)$.

\subsubsection{HS relation}\label{HS} This corresponds to Proposition~\ref{HSrelation}.
\begin{multline*}
	\begin{tikzpicture}
		\draw [bgd] (0,0) circle [radius=1.0];
		\draw [repV] (45:1) to [out=225,in=105] (225:0.5);
		\draw [blue!20, line width =3mm] (135:1) to [out=315,in=75] (315:0.5);
		\draw [repV] (135:1) to [out=315,in=75] (315:0.5);
		\draw [repW] (225:0.5) to [out=345,in=195] (315:0.5);
		\draw [repV] (225:0.5) -- (225:1);
		\draw [repV] (315:0.5) -- (315:1);
	\end{tikzpicture}
	\raisebox{1cm}{\quad $= t^{-2}u^{-2}v^{-2}$\quad}
	\begin{tikzpicture}
		\draw [bgd] (0,0) circle [radius=1.0];
		\draw [repV] (45:1) to [out=225,in=135] (315:1);
		\draw [repV] (135:1) to [out=315,in=45] (225:1);
	\end{tikzpicture}
	\raisebox{1cm}{\quad $+t^2u^2v^2$\quad}
	\begin{tikzpicture}
		\draw [bgd] (0,0) circle [radius=1.0];
		\draw [repV] (45:1) to [out=225,in=315] (135:1);
		\draw [repV] (315:1) to [out=135,in=45] (225:1);
	\end{tikzpicture}  \\
	\raisebox{1cm}{\quad $- t^{-1}u^{-1}$\quad}
	\begin{tikzpicture}
		\draw [bgd] (0,0) circle [radius=1.0];
		\draw [repV] (45:1) to [out=225,in=30] (0,0.3);
		\draw [repV] (135:1) to [out=315,in=150] (0,0.3);
		\draw [repV] (225:1) to [out=45,in=210] (0,-0.3);
		\draw [repV] (315:1) to [out=135,in=330] (0,-0.3);
		\draw [repW] (0,0.3) -- (0,-0.3);
	\end{tikzpicture}
	\raisebox{1cm}{\quad $-tu$\quad}
	\begin{tikzpicture}
		\draw [bgd] (0,0) circle [radius=1.0];
		\draw [repV] (45:1) to [out=225,in=60] (0.3,0);
		\draw [repV] (135:1) to [out=315,in=120] (-0.3,0);
		\draw [repV] (225:1) to [out=45,in=240] (-0.3,0);
		\draw [repV] (315:1) to [out=135,in=300] (0.3,0);
		\draw [repW] (0.3,0) -- (-0.3,0);
	\end{tikzpicture}
	\raisebox{1cm}{\quad $+$\quad}
	\begin{tikzpicture}
		\draw [bgd] (0,0) circle [radius=1.0];
		\draw [repV] (45:1) -- (225:1);
		\draw [blue!20, line width=3mm] (135:1) -- (315:1);
		\draw [repV] (135:1) -- (315:1);
	\end{tikzpicture}
\end{multline*}

\subsubsection{Triangle relation}
The triangle relation \eqref{triangle} is:
\begin{equation}
	\begin{tikzpicture}
		\draw [blue, fill=blue!20] (0,0) circle [radius=1.0];
		\draw [cyan, ultra thick] (90:0.5) -- (90:1);
		\draw [red, ultra thick] (210:0.5) -- (210:1);
		\draw [red, ultra thick] (330:0.5) -- (330:1);
		\draw [red, ultra thick] (90:0.5) to [out=210, in=90] (210:0.5);
		\draw [cyan, ultra thick] (210:0.5) to [out=330, in=210] (330:0.5);
		\draw [red, ultra thick] (330:0.5) to [out=90, in=330] (90:0.5);
	\end{tikzpicture}
	\raisebox{1cm}{\quad $=\frac{[1,1,2]}{[2,2,4]}\left(\frac{[4,4,0]}{[2,2,0]}-\psi\right)$}
	\begin{tikzpicture}
		\draw [blue, fill=blue!20] (0,0) circle [radius=1.0];
		\draw [cyan, ultra thick] (0,0) -- (90:1);
		\draw [red, ultra thick] (0,0) -- (210:1);
		\draw [red, ultra thick] (0,0) -- (330:1);
	\end{tikzpicture}
\end{equation}
This uses
\begin{equation}\label{eq:tau}
	\alpha^{-2}\gamma+\alpha^2\gamma^{-1} = \frac{[4,4,0]}{[2,2,0]}
	\qquad
	\alpha^{-2}\gamma^2+\alpha^2\gamma^{-2} = \frac{[2,2,4]}{[1,1,2]}
\end{equation}

\subsubsection{Skein relation}\label{skein} This corresponds to Proposition~\ref{skeinrelation}.
\begin{multline*}
	\begin{tikzpicture}[baseline=0]
		\draw [bgd] (0,0) circle [radius=1.0];
		\draw [repV] (45:1) -- (225:1);
		\draw [blue!20, line width=3mm] (135:1) -- (315:1);
		\draw [repV] (135:1) -- (315:1);
	\end{tikzpicture}
	\parbox{1cm}{\quad $-$\quad}
	\begin{tikzpicture}[baseline=0]
		\draw [bgd] (0,0) circle [radius=1.0];
		\draw [repV] (135:1) -- (315:1);
		\draw [blue!20, line width=3mm] (45:1) -- (225:1);
		\draw [repV] (45:1) -- (225:1);
	\end{tikzpicture}
	\parbox{5cm}{\quad $=\frac{[1,0,0][2,4,2]}{[3,3,2][1,2,1]}$} \\
	\left[\parbox{2cm}{$[1,1,1]$}
	\left(
	\begin{tikzpicture}[baseline=0]
		\draw [bgd] (0,0) circle [radius=1.0];
		\draw [repV] (45:1) to [out=225,in=135] (315:1);
		\draw [repV] (135:1) to [out=315,in=45] (225:1);
	\end{tikzpicture}
	\parbox{1cm}{\quad $-$\quad}
	\begin{tikzpicture}[baseline=0]
		\draw [bgd] (0,0) circle [radius=1.0];
		\draw [repV] (45:1) to [out=225,in=315] (135:1);
		\draw [repV] (315:1) to [out=135,in=45] (225:1);
	\end{tikzpicture}
	\right) 
	\parbox{2cm}{\quad $-[2,2,1]$\quad}
	\left(
	\begin{tikzpicture}[baseline=0]
		\draw [bgd] (0,0) circle [radius=1.0];
		\draw [repV] (45:1) to [out=225,in=30] (0,0.3);
		\draw [repV] (135:1) to [out=315,in=150] (0,0.3);
		\draw [repV] (225:1) to [out=45,in=210] (0,-0.3);
		\draw [repV] (315:1) to [out=135,in=330] (0,-0.3);
		\draw [repW] (0,0.3) -- (0,-0.3);
	\end{tikzpicture}
	\parbox{1cm}{\quad $-$\quad}
	\begin{tikzpicture}[baseline=0]
		\draw [bgd] (0,0) circle [radius=1.0];
		\draw [repV] (45:1) to [out=225,in=60] (0.3,0);
		\draw [repV] (135:1) to [out=315,in=120] (-0.3,0);
		\draw [repV] (225:1) to [out=45,in=240] (-0.3,0);
		\draw [repV] (315:1) to [out=135,in=300] (0.3,0);
		\draw [repW] (0.3,0) -- (-0.3,0);
	\end{tikzpicture}
	\right)\right]
\end{multline*}

\subsubsection{Square relation}\label{square} This corresponds to Proposition~\ref{squarerelation}.

\begin{multline*}
	\begin{tikzpicture}[baseline=0]
		\draw [bgd] (0,0) circle [radius=1.0];
		\draw [repV] (45:0.5) -- (45:1);
		\draw [repV] (135:0.5) -- (135:1);
		\draw [repV] (225:0.5) -- (225:1);
		\draw [repV] (315:0.5) -- (315:1);
		\draw [repW] (45:0.5) to [out=165,in=15] (135:0.5);
		\draw [repV] (135:0.5) to [out=255,in=105] (225:0.5);
		\draw [repW] (225:0.5) to[out=345,in=195] (315:0.5);
		\draw [repV] (315:0.5) to [out=75,in=285] (45:0.5);
	\end{tikzpicture}
	\parbox{1.8cm}{\quad $= Sq_I$\quad}
	\begin{tikzpicture}[baseline=0]
		\draw [bgd] (0,0) circle [radius=1.0];
		\draw [repV] (45:1) to [out=225,in=135] (315:1);
		\draw [repV] (135:1) to [out=315,in=45] (225:1);
	\end{tikzpicture}
	\parbox{1.8cm}{\quad $+ Sq_U$\quad}
	\begin{tikzpicture}[baseline=0]
		\draw [bgd] (0,0) circle [radius=1.0];
		\draw [repV] (45:1) to [out=225,in=315] (135:1);
		\draw [repV] (315:1) to [out=135,in=45] (225:1);
	\end{tikzpicture} \\
	\parbox{1.8cm}{\quad $+ Sq_K$\quad}
	\begin{tikzpicture}[baseline=0]
		\draw [bgd] (0,0) circle [radius=1.0];
		\draw [repV] (45:1) to [out=225,in=30] (0,0.3);
		\draw [repV] (135:1) to [out=315,in=150] (0,0.3);
		\draw [repV] (225:1) to [out=45,in=210] (0,-0.3);
		\draw [repV] (315:1) to [out=135,in=330] (0,-0.3);
		\draw [repW] (0,0.3) -- (0,-0.3);
	\end{tikzpicture}
	\parbox{1.8cm}{\quad $+ Sq_H$\quad}
	\begin{tikzpicture}[baseline=0]
		\draw [bgd] (0,0) circle [radius=1.0];
		\draw [bgd] (0,0) circle [radius=1.0];
		\draw [repV] (45:1) to [out=225,in=60] (0.3,0);
		\draw [repV] (135:1) to [out=315,in=120] (-0.3,0);
		\draw [repV] (225:1) to [out=45,in=240] (-0.3,0);
		\draw [repV] (315:1) to [out=135,in=300] (0.3,0);
		\draw [repW] (0.3,0) -- (-0.3,0);
	\end{tikzpicture}
	\parbox{1.8cm}{\quad $+Sq_S$\quad}
	\begin{tikzpicture}[baseline=0]
		\draw [bgd] (0,0) circle [radius=1.0];
		\draw [repV] (45:1) -- (225:1);
		\draw [blue!20, line width=3mm] (135:1) -- (315:1);
		\draw [repV] (135:1) -- (315:1);
	\end{tikzpicture}
\end{multline*}
where
\begin{equation*}
	Sq_S = \frac{[3,3,2]^2}%
	{ [1,0,0][2,2,1]}\:
	\frac{[1,2,1]}{[2,4,2]}
\end{equation*}

\subsubsection{Bigon relation}\label{digon} This corresponds to Proposition~\ref{digonrelation}.

\begin{multline*}
	\begin{tikzpicture}
		\draw [bgd] (0,0) circle [radius=1.0];
		\draw [repV] (45:1) to [out=225,in=90] (-0.5,0) to [out=270,in=135] (315:1);
		\draw [blue!20, line width =3mm] (135:1) to [out=315,in=90] (0.5,0) to [out=270,in=45] (225:1);
		\draw [repV] (135:1) to [out=315,in=90] (0.5,0) to [out=270,in=45] (225:1);
	\end{tikzpicture}
	\raisebox{1cm}{\quad $=1+u\,t\,\frac{[1,0,0][2,2,1]}{[3,3,2]}
		\frac{[2,4,2]}{[1,2,1]}$\quad}
	\begin{tikzpicture}
		\draw [bgd] (0,0) circle [radius=1.0];
		\draw [repV] (45:1) to [out=225,in=135] (315:1);
		\draw [repV] (135:1) to [out=315,in=45] (225:1);
	\end{tikzpicture} \\
	\raisebox{1cm}{\quad $-\,t^5\,u^5v^5\,\frac{[1,0,0][2,4,2]}{[3,3,3][1,2,1]}
		\left(\frac{([6,6,6][1,1,1])}{[2,2,2][3,3,3]}-t^2\,u^2\right)$\quad}
	\begin{tikzpicture}
		\draw [bgd] (0,0) circle [radius=1.0];
		\draw [repV] (45:1) to [out=225,in=315] (135:1);
		\draw [repV] (315:1) to [out=135,in=45] (225:1);
	\end{tikzpicture}  \\
	\raisebox{1cm}{\quad $-t^4\,u^4\,v^2\,\frac{[1,0,0][2,2,1]}{[3,3,2]}
		\frac{[2,4,2]}{[1,2,1]}$\quad}
	\begin{tikzpicture}
		\draw [bgd] (0,0) circle [radius=1.0];
		\draw [repV] (45:1) to [out=225,in=30] (0,0.3);
		\draw [repV] (135:1) to [out=315,in=150] (0,0.3);
		\draw [repV] (225:1) to [out=45,in=210] (0,-0.3);
		\draw [repV] (315:1) to [out=135,in=330] (0,-0.3);
		\draw [repW] (0,0.3) -- (0,-0.3);
	\end{tikzpicture} \\
	\raisebox{1cm}{\quad $-t\,u\,v^2\,\frac{[1,0,0][2,2,1]}{[3,3,2]}
		\frac{[2,2,0]]}{[1,1,0]}\frac{[2,4,2]}{[1,2,1]}$\quad}
	\begin{tikzpicture}
		\draw [bgd] (0,0) circle [radius=1.0];
		\draw [repV] (45:1) to [out=225,in=60] (0.3,0);
		\draw [repV] (135:1) to [out=315,in=120] (-0.3,0);
		\draw [repV] (225:1) to [out=45,in=240] (-0.3,0);
		\draw [repV] (315:1) to [out=135,in=300] (0.3,0);
		\draw [repW] (0.3,0) -- (-0.3,0);
	\end{tikzpicture} \\
	\raisebox{1cm}{\quad $+v\,u^2\,t^2\,[1,0,0]\frac{[2,4,2]}{[1,2,1]}$\quad}
	\begin{tikzpicture}
		\draw [bgd] (0,0) circle [radius=1.0];
		\draw [repV] (45:1) -- (225:1);
		\draw [blue!20, line width=3mm] (135:1) -- (315:1);
		\draw [repV] (135:1) -- (315:1);
	\end{tikzpicture}
\end{multline*}

This gives a linear ribbon category which satisfies the Assumption~\ref{assumption}.

The category $\Dcat$ has two involutions. The \Dfn{bar involution} is anti-linear
for the bar involution on $D(t,u,v)$ and is given on diagrams by reversing all crossings.
The second involution is anti-linear for the involution on $D(t,u,v)$ given by
$t\mapsto t$, $u\mapsto v$, $v\mapsto u$. It is defined on objects by
$V\mapsto V$, $Y\mapsto W$, $W\mapsto Y$, $L\mapsto L$, $X\mapsto X$.
This follows from the universal property.

These two involutions commute.

\subsection{Calculations}\label{sec:2point}

The structure constants of the algebra $A(2)$ with respect to the basis $\bB$,
defined in \eqref{basis}, have been given in \S~\ref{algebra}. This algebra
is commutative by construction and can be checked to be associative by a
direct calculation.

The algebra $A(2)$ is generically semisimple with a basis of orthogonal
idempotents. We introduce the images of these idempotents as additional objects
and write
\begin{equation*}
	V\otimes V \cong 1 \oplus Y \oplus W \oplus L \oplus X 
\end{equation*}
These idempotents can be calculated and this gives the eigenvalues
of $S$. These eigenvalues are shown in Figure~\ref{braid:generic}.

\begin{figure}
	\begin{tabular}{ccccc}
		$1$ & $L$ & $Y$ & $X$ & $W$ \\ \hline
		$t^6\,u^6\,v^6$ & $-1$ & $u^2$ & $-t^4\,u^2\,v^2$ & $v^2$ 
	\end{tabular}
	\caption{Eigenvalues}\label{braid:generic}
\end{figure}

The quantum dimension of the image of an idempotent $\pi$ is
calculated as $\dim_q(\pi) = \delta_V\,\Tr(\pi)$ where $\Tr$ is the trace map
defined in Definition~\ref{defn:tr}.

The quantum dimensions are given by:
\begin{align*}
\dim_q (L) &=
\frac{[1,2,2][2,2,1][2,3,3][3,3,2][2,3,2][3,3,3] [2,0,2][1,2,1)^2[2,2,4][4,2,2][4,2,4]}%
{[1,2,1][2,0,0][0,0,2][1,0,0][0,0,1][1,1,1]
[1,0,1][2,4,2)^2[1,1,2][2,1,1][2,1,2]}	\\
\dim_q(Y)	&=
\frac{[0,3,3][1,2,2][0,2,1][2,2,2][3,3,3][2,3,2]
[1,2,1][4,2,4][2,2,4][3,3,2]}%
{[2,0,0][1,0,0][-1,0,1][0,0,1][0,1,1][0,4,2]
[1,1,1][2,4,2][2,1,2][1,1,2]} \\
\dim_q(X)	&= (-1)\frac{[3,3,2][2,3,3]}{[1,0,0][0,0,1]} \frac{[0,6,0][0,1,0]}{[0,2,0][0,3,0]}\\
&\qquad\qquad\times \frac{[1,2,1]^2[2,2,2][3,3,3][2,0,2][1,2,0][0,2,1][2,2,4][4,2,2]}%
{[2,4,2]^2[1,1,1]^2[1,0,1][2,4,0][0,4,2][1,1,2][2,1,1]}
\end{align*}

\subsection{Classical category}\label{sec:classical}
In this section we specialise using $\ev$ and construct a rigid symmetric monoidal
category. This is a quotient of the category $\bQ^{(0)}(m,n,p)\Gcat$,
the free $\bQ^{(0)}(m,n,p)$-linear category on the category $\Gcat$, discussed in \S~\ref{sec:Gcat}.

For the coefficients which have been written as elements of $D(t,u,v)$ it is
straightforward to apply $\ev$. The classical version of the digon relation
is $S^2=1$ and it is clear that this is obtained by applying $\ev$. The only relation
that requires discussion is the square relation.

The relations in \S~\ref{initial} hold with
\begin{equation*}
	\delta_V = (-2)\frac{ (2m + 3n + 3p)(3m + 3n + 2p)}{pm}
\end{equation*}
\begin{equation*}
	\delta_W = \frac{9(m+n+p)(2m+2n+p)(2m+3n+2p)(2m+3n+3p)}%
	{mp^2(m-p)}
\end{equation*}

\subsubsection{HS relation}\label{HSc}
The classical version of the $HS$  relation in \S~\ref{HS} is
\begin{multline*}
	\begin{tikzpicture}
		\draw [bgd] (0,0) circle [radius=1.0];
		\draw [repV] (45:1) to [out=225,in=105] (225:0.5);
		\draw [repV] (135:1) to [out=315,in=75] (315:0.5);
		\draw [repW] (225:0.5) to [out=345,in=195] (315:0.5);
		\draw [repV] (225:0.5) -- (225:1);
		\draw [repV] (315:0.5) -- (315:1);
	\end{tikzpicture}
	\raisebox{1cm}{\quad $= $\quad}
	\begin{tikzpicture}
		\draw [bgd] (0,0) circle [radius=1.0];
		\draw [repV] (45:1) to [out=225,in=135] (315:1);
		\draw [repV] (135:1) to [out=315,in=45] (225:1);
	\end{tikzpicture}
	\raisebox{1cm}{\quad $+$\quad}
	\begin{tikzpicture}
		\draw [bgd] (0,0) circle [radius=1.0];
		\draw [repV] (45:1) to [out=225,in=315] (135:1);
		\draw [repV] (315:1) to [out=135,in=45] (225:1);
	\end{tikzpicture}  \\
	\raisebox{1cm}{\quad $- $\quad}
	\begin{tikzpicture}
		\draw [bgd] (0,0) circle [radius=1.0];
		\draw [repV] (45:1) to [out=225,in=30] (0,0.3);
		\draw [repV] (135:1) to [out=315,in=150] (0,0.3);
		\draw [repV] (225:1) to [out=45,in=210] (0,-0.3);
		\draw [repV] (315:1) to [out=135,in=330] (0,-0.3);
		\draw [repW] (0,0.3) -- (0,-0.3);
	\end{tikzpicture}
	\raisebox{1cm}{\quad $-$\quad}
	\begin{tikzpicture}
		\draw [bgd] (0,0) circle [radius=1.0];
		\draw [repV] (45:1) to [out=225,in=60] (0.3,0);
		\draw [repV] (135:1) to [out=315,in=120] (-0.3,0);
		\draw [repV] (225:1) to [out=45,in=240] (-0.3,0);
		\draw [repV] (315:1) to [out=135,in=300] (0.3,0);
		\draw [repW] (0.3,0) -- (-0.3,0);
	\end{tikzpicture}
	\raisebox{1cm}{\quad $+$\quad}
	\begin{tikzpicture}
		\draw [bgd] (0,0) circle [radius=1.0];
		\draw [repV] (45:1) -- (225:1);
		\draw [repV] (135:1) -- (315:1);
	\end{tikzpicture}
\end{multline*}

\subsubsection{Bigons}
Define $\phi,\psi$ by
\begin{equation*}
	\phi = (-3)\frac{(m + n + p)(2m + 3n + 2p)}%
	{mp}
\end{equation*}
\begin{equation*}
	\psi = \frac23 \frac{(m - p)(3m + 3n + 2p)}%
	{m(2m + 2n + p)}
\end{equation*}
\begin{equation}
	\begin{tikzpicture}[baseline=0]
		\draw [bgd] (0,0) circle [radius=1.0];
		\draw [repW] (0,0.5) to [out=210, in=90] (-0.4,0) to [out=270, in=150] (0,-0.5);
		\draw [repV] (0,0.5) to [out=330, in=90] (0.4,0) to [out=270, in=30](0,-0.5);
		\draw [repV] (0,0.5) -- (0,1);
		\draw [repV] (0,-0.5) -- (0,-1);
	\end{tikzpicture}
	\parbox[b]{5cm}{\quad$=(-3)\frac{(m + n + p)(2m + 3n + 2p)}{mp}$}
	\begin{tikzpicture}[baseline=0]
		\draw [bgd] (0,0) circle [radius=1.0];
		\draw [repV] (0,1) -- (0,-1);
	\end{tikzpicture}
\end{equation}

\begin{equation}
	\begin{tikzpicture}[baseline=0]
		\draw [bgd] (0,0) circle [radius=1.0];
		\draw [repV] (0,0.5) to [out=210, in=90] (-0.4,0) to [out=270, in=150] (0,-0.5);
		\draw [repV] (0,0.5) to [out=330, in=90] (0.4,0) to [out=270, in=30](0,-0.5);
		\draw [repW] (0,0.5) -- (0,1);
		\draw [repW] (0,-0.5) -- (0,-1);
	\end{tikzpicture}
	\parbox[b]{5cm}{\quad$= \frac23 \frac{(m - p)(3m + 3n + 2p)}{m(2m + 2n + p)}$}
	\begin{tikzpicture}[baseline=0]
		\draw [bgd] (0,0) circle [radius=1.0];
		\draw [repW] (0,1) -- (0,-1);
	\end{tikzpicture}
\end{equation}

\subsubsection{Triangle relation}
The coefficient defined in \eqref{triangle} is
\begin{equation*}
	\tau = 
	\frac{(3m^2 + 3mn + 4mp + 3np + 2p^2)}{3m(2m + 2n + p)}
\end{equation*}

\subsubsection{Bigon relation}
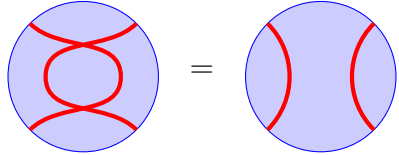
\begin{equation*}
	\begin{tikzpicture}
		\draw [bgd] (0,0) circle [radius=1.0];
		\draw [repV] (45:1) to [out=225,in=90] (-0.5,0) to [out=270,in=135] (315:1);
		\draw [repV] (135:1) to [out=315,in=90] (0.5,0) to [out=270,in=45] (225:1);
	\end{tikzpicture}
	\raisebox{1cm}{\quad $=$\quad}
	\begin{tikzpicture}
		\draw [bgd] (0,0) circle [radius=1.0];
		\draw [repV] (45:1) to [out=225,in=135] (315:1);
		\draw [repV] (135:1) to [out=315,in=45] (225:1);
	\end{tikzpicture} 
\end{equation*}

\subsubsection{Square relation}\label{rel:square}
Since we have not expressed the coefficients in the square relation as
elements of $D(t,u,v)$ we cannot ealuate directly. Instead
we evaluate $Sq_S$ directly and obtain the remaining coefficients
by solving a linear system of equations. These equations are
the associativity conditions
\begin{equation*}
	(HH)U = H(HU) \qquad (HH)K = H(HK) \qquad (HH)S = H(HS)
\end{equation*}
together with the equation, which is given by applying the trace map to the relation,
\begin{equation*}
 \psi\phi = Sq_I\delta_V + Sq_U + Sq_K\phi + Sq_S
\end{equation*}

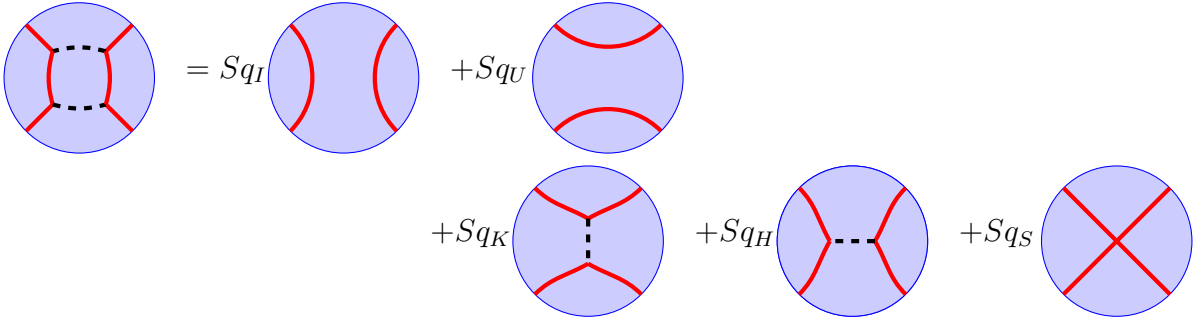
\begin{multline*}
	\begin{tikzpicture}[baseline=0]
		\draw [bgd] (0,0) circle [radius=1.0];
		\draw [repV] (45:0.5) -- (45:1);
		\draw [repV] (135:0.5) -- (135:1);
		\draw [repV] (225:0.5) -- (225:1);
		\draw [repV] (315:0.5) -- (315:1);
		\draw [repW] (45:0.5) to [out=165,in=15] (135:0.5);
		\draw [repV] (135:0.5) to [out=255,in=105] (225:0.5);
		\draw [repW] (225:0.5) to[out=345,in=195] (315:0.5);
		\draw [repV] (315:0.5) to [out=75,in=285] (45:0.5);
	\end{tikzpicture}
	\parbox{1.5cm}{\quad $= Sq_I$\quad}
	\begin{tikzpicture}[baseline=0]
		\draw [bgd] (0,0) circle [radius=1.0];
		\draw [repV] (45:1) to [out=225,in=135] (315:1);
		\draw [repV] (135:1) to [out=315,in=45] (225:1);
	\end{tikzpicture}
	\parbox{1.5cm}{\quad $+ Sq_U$\quad}
	\begin{tikzpicture}[baseline=0]
		\draw [bgd] (0,0) circle [radius=1.0];
		\draw [repV] (45:1) to [out=225,in=315] (135:1);
		\draw [repV] (315:1) to [out=135,in=45] (225:1);
	\end{tikzpicture} \\
	\parbox{1.5cm}{\quad $+ Sq_K$\quad}
	\begin{tikzpicture}[baseline=0]
		\draw [bgd] (0,0) circle [radius=1.0];
		\draw [repV] (45:1) to [out=225,in=30] (0,0.3);
		\draw [repV] (135:1) to [out=315,in=150] (0,0.3);
		\draw [repV] (225:1) to [out=45,in=210] (0,-0.3);
		\draw [repV] (315:1) to [out=135,in=330] (0,-0.3);
		\draw [repW] (0,0.3) -- (0,-0.3);
	\end{tikzpicture}
	\parbox{1.5cm}{\quad $+ Sq_H$\quad}
	\begin{tikzpicture}[baseline=0]
		\draw [bgd] (0,0) circle [radius=1.0];
		\draw [bgd] (0,0) circle [radius=1.0];
		\draw [repV] (45:1) to [out=225,in=60] (0.3,0);
		\draw [repV] (135:1) to [out=315,in=120] (-0.3,0);
		\draw [repV] (225:1) to [out=45,in=240] (-0.3,0);
		\draw [repV] (315:1) to [out=135,in=300] (0.3,0);
		\draw [repW] (0.3,0) -- (-0.3,0);
	\end{tikzpicture}
	\parbox{1.5cm}{\quad $+Sq_S$\quad}
	\begin{tikzpicture}[baseline=0]
		\draw [bgd] (0,0) circle [radius=1.0];
		\draw [repV] (45:1) -- (225:1);
		\draw [repV] (135:1) -- (315:1);
	\end{tikzpicture}
\end{multline*}
where
\begin{gather*}
	Sq_I =  \frac{(5m^2 + 8mn + 3n^2 + 4mp + 2np)}%
	{2m(2m + 2n + p)} \\
	Sq_U =  (-1)\frac{(4m^2 + 10mn + 6n^2 + 5mp + 7np + 2p^2)(3m + 3n + 2p)}%
	{2mp(2m + 2n + p)} \\
	Sq_K = \frac{(4m + 3n + 2p)(3m + 3n + 2p)}{3m(2m + 2n + p)} \\
	Sq_H = (-1)\frac{(5m^2 + 11mn + 6n^2 + 7mp + 7np + 2p^2)}%
	{m(2m + 2n + p)} \\
	Sq_S =\frac{(3m + 3n + 2p)^2}{2m(2m + 2n + p)}
\end{gather*}

\subsubsection{Dimensions}
The dimensions are
\begin{equation*}
	\dim Y = 
\frac{(m + n + p)(m + 2n + 2p)(2m + 3n + 2p)(3m + 3n + 2p}%
{m^2p(-m + p)}
\end{equation*}

\begin{equation*}
\dim X = 3\frac{(m + 2n + 2p)(2m + 2n + p)(2m + 3n + 2p)(2m + 3n + 3p)(3m + 3n + 2p)}%
{m^2p^2 (m + 2n + p) }
\end{equation*}

\begin{equation*}
\dim L = 3\frac{(m + 2n + 2p)(2m + 2n + p)(2m + 3n + 2p)(2m + 3n + 3p)(3m + 3n + 2p)}%
{p^2 m^2 (m + 2n + p)}
\end{equation*}

\section{Evaluation functors}\label{sec:functors}
In this section we discuss the evaluation functors. These are organised into lines.
For each line we have a homomorphism $D(t,u,v)\to D(q,z)$ and this gives a 
$D(q,z)$-linear ribbon category by taking the specialisation $\Dcat$

\subsection{Quantum groups}
Let $C$ be a finite type Cartan matrix and let $\fg(C)$ be the associated
semisimple complex Lie algebra. Let $\Mcat_\bC(C)$ be the category of finite dimensional
$\fg(C)$-modules. This is a rigid symmetric monoidal category which we regard as a
ribbon category. In this section we follow \cite[\S~5]{Savage2022} and construct a flat deformation of this category as a ribbon category. Let $A$ be the local ring defined in \S~\ref{scalars}. The homomorphism $A\to \bC$ is $A\to \bQ \to\bC$.
The flat deformation is a ribbon category, $\Mcat_A(C)$, with the properties
\begin{itemize}
	\item The category $\Mcat_A(C)$ is enriched in the category of finitely generated free $A$-modules
	\item The category $\Mcat_A(C)$ is semisimple abelian.
	\item The category $\Mcat_A(C)\otimes_A\bC$ is equivalent to $\Mcat_\bC(C)$.
	\item The equivalence $\Mcat_A(C)\otimes_A\bC\to \Mcat_\bC(C)$ is an
	equivalence of $\bC$-linear ribbon categories.
\end{itemize}

This property is apparent in the Drinfel'd approach to quantum groups,
 \cite{Drinfeld1987}. by design.
In this approach the quantum group is a ribbon Hopf algebra in the category of
complete $\bC[[h]]$-modules. Here $\bC[[h]]$ is
the ring of formal power series and complete means complete in the $h$-adic topology.
Here we take the algebraic approach to quantum groups presented in \cite{Lusztig1993}, \cite{Chari1995}, \cite{Jantzen1996}.

The isomorphism classes of simple $\fg(C)$-modules are parametrised by dominant weights.
Let $L_\bC(\lambda)$ be the simple module associated to the dominant weight $\lambda$.
For each $\lambda$ there is a simple object $L_A(\lambda)$ which is characterised by the properties
\begin{itemize}
	\item $L_A(\lambda)$ is finitely generated and free as an $A$-module.
	\item $L_A(\lambda)\otimes_A\bC\cong L_\bC(\lambda)$.
\end{itemize}

The twist for the Drinfeld quantum group is implicitly given in \cite[(2.22)]{Reshetikhin1987}. A discussion for $U_q(C)$ is given in \cite[\S~2]{Leduc1997}.

The twist $\theta_{L_A(\lambda)}$ is defined as an element of $\End(L_A(\lambda))$.
Since we have a canonical identification $\End(L_A(\lambda))\cong A$ we consider
$\theta_{L_A(\lambda)}$ as an element of $A$.
\begin{prop} The twist $\theta_{L_A(\lambda)}\in A$ is
	\begin{equation*}
\theta_{L_A(\lambda)} = q^{\left\langle \lambda,\lambda+\rho\right\rangle}		
	\end{equation*}
\end{prop}

Finally, we evaluate the quantum dimensions using the quantum Weyl dimension
formula. This is folk-lore and an outline proof is \cite[Lemma~5.7]{Savage2022}.
\begin{lemma}[Quantum Weyl dimension formula]
	We have
	\begin{equation} \label{qWeyl}
		\dim_q L_A(\lambda)
		= \prod_{\nu \in \Phi^+} \frac{[(\lambda + \rho, \nu)]}{[(\rho,\nu)]}.
	\end{equation}
\end{lemma} 

This can be extended to allow diagram automorphisms. Recall that a \Dfn{diagram automorphism},
$g$, is a permutation of the vertices of the Dynkin diagram such that $a_{g(i),g(j)}= a_{i,j}$
for all vertices $i$ and $j$. The set of diagram automorphisms is a finite group and we let
$\Gamma$ be a subgroup of this group. This group acts on the enveloping algebra
$U(\fg(C))$ and so we can form the semidirect  product with the group algebra,
$U(\fg(C))\rtimes \bC\Gamma$.  Let $\Mcat_\bC(C,\Gamma)$ be the category of finite dimensional
$U(\fg(C))\rtimes \bC\Gamma$-modules. This is a $\bC$-linear rigid symmetric monoidal category which also has a flat deformation as a ribbon category. This is a ribbon category, $\Mcat_A(C)$,
with the properties 
\begin{itemize}
	\item The category $\Mcat_A(C,\Gamma)$ is enriched in the category of finitely generated free $A$-modules
	\item The category $\Mcat_A(C,\Gamma)$ is semisimple abelian.
	\item The category $\Mcat_A(C)\otimes_A\bC$ is equivalent to $\Mcat_\bC(C,\Gamma)$.
	\item The equivalence $\Mcat_A(C,\Gamma)\otimes_A\bC\to \Mcat_\bC(C,\Gamma)$ is an
	equivalence of $\bC$-linear ribbon categories.
\end{itemize}

There is a restriction functor $\Mcat_A(C,\Gamma)\to \Mcat_A(C)$ which we denote by
$M\mapsto M\downarrow$.

The group $\Gamma$ acts on dominant weights. The isomorphism classes of simple $\fg(C)$-modules are parametrised by the orbits of the action of $\Gamma$ on dominant weights.
Let $L_A(\mathcal{O})$ be the simple module associated to the orbit $\mathcal{O}$.
The restriction of this module is given by
\begin{equation*}
	L_A(\mathcal{O})\downarrow \cong \bigoplus_{\lambda\in\mathcal{O}} L_A(\lambda)
\end{equation*}

The quantum dimensions satisfy $\dim_q M = \dim_q M\downarrow$. The function
$\lambda\mapsto \dim_q L_A(\lambda)$ is constant on each orbit so we have
\begin{equation*}
		\dim_q L_A(\mathcal{O}) = |\mathcal{O}|\,\dim_q L_A(\lambda) 
		\qquad\text{for any $\lambda\in\mathcal{O}$.}
\end{equation*}

The function
$\lambda\mapsto \theta_{L_A(\lambda)}$ is constant on each orbit so we have
\begin{equation*}
	\theta_{L_A(\mathcal{O})} = \theta_{L_A(\lambda)}
	\qquad\text{for any $\lambda\in\mathcal{O}$.}
\end{equation*}

Let $\lambda$ be a dominant weight fixed by $\Gamma$ so $\mathcal{O}=\{\lambda\}$
and $L_A(\mathcal{O})\downarrow \cong L_A(\lambda)$. There are examples where
$L_A(\mathcal{O})\in \Mcat_A(C,\Gamma)$ satisfies the Assumption~\ref{assumption} but
$L_A(\lambda)\in \Mcat_A(C)$ does not.

\subsection{Evaluation functors}
For some Cartan matrix, $C$, let $\lambda$ be a dominant weight such that
the symmetric and exteroior powers of $L_\bC(\lambda)$ have the following decompositions
\begin{equation}\label{condition}
	S^2 L_\bC(\lambda) \cong L_\bC(0) \oplus L_\bC(2\lambda) \oplus L_\bC(\mu)
	\qquad
	\bigwedge\nolimits^2 L_\bC(\lambda) \cong L_\bC(\theta) \oplus L_\bC(\nu)
\end{equation}
Here $0$ is the zero weight so $L(0)$ is the trivial representation and $\theta$ is the highest root so $L(\theta)$ is the adjoint representation.

For any $\lambda$, $L_\bC(2\lambda)$ is a summand of $S^2 L_\bC(\lambda)$.
The condition that $L_\bC(0)$ is a summand of $S^2 L_\bC(\lambda)$
is equivalent to the condition that $L_\bC(\lambda)$ is self-dual with a 
non-degenerate symmetric inner product. In this case the adjoint representation is a summand
of $\bigwedge\nolimits^2 L_\bC(\lambda)$. Hence the condition \eqref{condition}
is equivalent to the conditions
\begin{itemize}
	\item $L_\bC(\lambda)$ is self-dual with a non-degenerate symmetric inner product
	\item the complement of $L_\bC(0) \oplus L_\bC(2\lambda) \subset S^2 L_\bC(\lambda)$
	is simple
	\item the complement of $L_\bC(\theta) \subset \bigwedge\nolimits^2 L_\bC(\lambda)$
is simple	
\end{itemize}

\begin{lemma} Let $L_\bC(\lambda)$ satisfy \eqref{condition}. Then
	\begin{itemize}
\item $V=L_A(\lambda)$ and $W=L_A(\mu)$ satisfy Assumption~\ref{assumption}
\item $V=L_A(\lambda)$ and $W=L_A(2\lambda)$ satisfy Assumption~\ref{assumption}
	\end{itemize}
\end{lemma}

\begin{proof}
Write $X$ and $\Rot_\bC$ for the actions of
 $S$ and $\Rot$ on the vector space $\End(L_\bC(\lambda)\otimes L_\bC(\lambda))$. Then $X^2=1$, by construction and $\Rot_A^2=1$ since
$L_\bC(\lambda)\otimes L_\bC(\lambda)$ is multiplicity free. These generate an action
of $\fS_3$ with $X=s_1$ and $\Rot_\bC=s_2s_1s_2$. Hence these are conjugate and the
$+1$ eigenspace of $\Rot_\bC$ has dimension 3 and the $-1$ eigenspace of $\Rot_\bC$ has dimension 2.

It remains to show that the set $\{1,U,K,H,S+S^{-1}\}$ is a basis of
$\End(L_A(\lambda)\otimes L_A(\lambda))$. It is sufficient to show that
the set $\{1,U,K,H,S+S^{-1}\}$ is a basis of $\End(L_\bC(\lambda)\otimes L_\bC(\lambda))$.
The proof follows the proof of \cite[Proposition~3.2]{Gandhi2022}.
\end{proof}

Let $L_\bC(\lambda)$ satisfy \eqref{condition}. Then put $V=L_A(\lambda)$
and $W=L_A(\mu)$. Then the twists are
\begin{equation*}
	\theta_V = q^{\left\langle \lambda,\lambda+\rho\right\rangle}
	\qquad	
	\theta_W = q^{\left\langle \mu,\mu+\rho\right\rangle}	
\end{equation*}

\begin{defn}
Let $\Tcat(t,u,v)$ be the category $\Tcat$ with parameters in $D(t,u,v)$,
where $\alpha,\beta,\gamma$ are given in \eqref{eq:twists},
$\delta_V$ is given in \eqref{eq:deltaV}, $\delta_W$ is given in \eqref{eq:deltaW},
$\psi$ is given in \eqref{eq:psi} and $\phi$ is given in \eqref{eq:phi}.
\end{defn}

\begin{thm} Let $\sigma\colon D(t,u,v)\to A$ be a homomorphism such that
\begin{align*}
	\sigma( tu) &=q^{\left\langle \mu,\mu+\rho\right\rangle / 2}
	&
	\sigma (v) &=  q^{\left\langle \lambda,\lambda+\rho\right\rangle}%
	q^{-\left\langle \mu,\mu+\rho\right\rangle / 2} \\
	\sigma(\delta_V) &= \dim_q L_A(\lambda)
	&
	\sigma(\delta_W) &= \dim_q L_A(\mu)
\end{align*}
Then there exists an $A$-linear ribbon functor
$\Dcat(t,u,v)\otimes_{D(t,u,v)} A\to \Mcat_A(C)$ with $V\mapsto L_A(\lambda)$
and $W\mapsto L_A(\mu)$. Furthermore the set of such functors has a free and transitive
action of $A^\times \times A^\times$.
\end{thm}
\begin{proof}
Then the free property gives an $A$-linear ribbon functor $\Tcat(t,u,v)\otimes_{D(t,u,v)} A\to \Mcat_A(C)$.
The results in \S then show that this factors through $\Dcat\otimes_{D(t,u,v)} A$
to give an $A$-linear ribbon functor $\Dcat(t,u,v)\otimes_{D(t,u,v)} A\to \Mcat_A(C)$.

Suppose $F$ is an $A$-linear ribbon functor $\Dcat(t,u,v)\otimes_{D(t,u,v)} A\to \Mcat_A(C)$ and $(c_V,c_W)$. Let $\mathrm{ev}_V$, $\mathrm{co}_V$, $\mathrm{ev}_W$,$\mathrm{co}_W$ be the morphisms making $V$ and $W$ self-dual. Define $(c_V,c_W)F$ by
\begin{align*}
(c_V,c_W)F \colon \mathrm{ev}_V &\mapsto c_V F(\mathrm{ev}_V)
&
(c_V,c_W)F \colon \mathrm{co}_V &\mapsto c_V^{-1} F(\mathrm{co}_V)
\\
(c_V,c_W)F \colon \mathrm{ev}_W &\mapsto c_W F(\mathrm{ev}_W)
&
(c_V,c_W)F \colon \mathrm{co}_W &\mapsto c_W^{-1} F(\mathrm{co}_W)
\end{align*}
and $(c_V,c_W)F = F$ on the remaining generators.

The relation \eqref{eq:bubbleV} will be satisfied with $\phi$ replaced by $\phi'$,
for some $\phi'$.
and the relation \eqref{eq:bubbleW} will be satisfied with $\psi$ replaced by $\psi'$,
for some $\psi'$.
Change $((c_V,c_W)F)(T)$ to $\frac{\phi}{\phi'}((c_V,c_W)F)(T)$.
Then both relations \eqref{eq:bubbleV} and \eqref{eq:bubbleW} are satisfied since
$\frac{\phi}{\phi'} = \frac{\psi}{\psi'}$.
\end{proof}
\subsection{Real series}
This series is the first line in the magic square associated to the real numbers.
The distinguishing feature is that $V\cong W$. This series has been studied in
\cite[Chapter~19]{Cvitanovic2008}, \cite[\S~8]{Landsberg2004a},\ 
The following is taken from \cite[unpublished]{Thurston2004}.
\begin{figure}[ht]
	\begin{tabular}{c|cccccccc}
		$k$ & -4/3 & -1 &-2/3 & 0 & 1 & 2 & 4 & 8 \\ \hline
		& $A_1$ & $OSp(1|2)$ & 1 & $\fS_3$ & $A_1$ & $A_2$ & $C_3$ & $F_4$ \\
	\end{tabular}
	\caption{Real line}
\end{figure}

The highest weights are shown in Figure~\ref{series:real}.
\begin{figure}[ht]
	\begin{tabular}{c|ccccc}
		& $V$ & $Y$ & $W$ & $L$ & $X$ \\ \hline
		$A_1$ & $4\omega_1$ & $8\omega_1$ & $4\omega_1$ & $2\omega_1$ & $6\omega_1$ \\
		$A_2\rtimes\fS_2$ & $\omega_1+\omega_2$ & $2\omega_1+2\omega_2$ & $\omega_1+\omega_2$ & $\omega_1+\omega_2$ & $3\omega_1^\ast$ \\
		$C_3$ & $\omega_2$ & $2\omega_2$ & $\omega_2$ & $2\omega_1$ & $\omega_1+\omega_3$ \\
		$F_4$ & $\omega_4$ & $2\omega_4$ & $\omega_4$ & $\omega_1$ & $\omega_3$ \\
	\end{tabular}
\caption{Highest weights}\label{series:real}
\end{figure}

This is the line $3m+3n+p=0$ and our choice of parametrisation is
$m=1$, $n=k-1$, $p=-3k$. The specialisation is
\begin{equation*}
	u \mapsto q^4 ,\quad t \mapsto z\,q^{-4} ,\quad v \mapsto z^{-3}
\end{equation*}

The values of the Casimir are given in Figure~\ref{casimir:real}.
Except that for $A_1$ there is a factor of 4.
\begin{figure}
	\begin{tabular}{ccccc}
		$V$ & $Y$ & $W$ & $L$ & $X$ \\ \hline
		$3k$ & $6k+4$ & $3k$ & $5k-4$ & $6k$
	\end{tabular}
	\caption{Casimir}\label{casimir:real}
\end{figure}

The eigenvalues of the braid matrix are given in Figure~\ref{braid:real}.
\begin{figure}
	\begin{tabular}{ccccc}
		$1$ & $L$ & $Y$ & $X$ & $W$ \\ \hline
		$z^{-12}$ & $-q^{-2}z^{-2}$ & $q^2$ & $-1$ & $z^{-6}$
	\end{tabular}
	\caption{Eigenvalues}\label{braid:real}
\end{figure}

The values of the parameter are
\begin{equation*}
	\alpha = t^{-4} \qquad \beta = t^{-6}  \qquad \gamma = t^{-3}
\end{equation*}

Then put $\delta_V = \delta = \delta_W$.
\begin{equation*}
	\delta = \frac{[-2k+10][4][k+1][k+6]}{[-k+5][2][2k+2]}
\end{equation*}

Then we also have $\phi=\psi$ and this is
\begin{equation*}
	\frac{[-2k+4][6][k+3][k+1]}{[2k+2][3][-k+2]}
\end{equation*}

The coefficient in the triangle relation is
\begin{equation*}
	\tau = \frac{[2]}{[4]}\left(\frac{[10]}{[5]}-\psi\right)
\end{equation*}

The coefficients in the skein relation are
\begin{equation*}
	z = (q-q^{-1})\frac{[2][2k+2]}{[3][k+1]} \qquad
	w = (q-q^{-1})\frac{[2k+2]}{[3][k+1]}
\end{equation*}

Taking the classical specialisation $\delta=12k+2$. This recovers the results
in \cite[Chapter~19]{Cvitanovic2008} with $n=12k+2$ and the results of \cite{Gandhi2022} with $\delta=12k+2$ and $\alpha=6k+2$.

Further specialisation by $t\mapsto q^{2}$, $[nk+r]\mapsto [n+2r]$ gives
\begin{equation*}
	\phi = \frac{[12][7]}{[6]} \qquad
	z = (q-q^{-1})\frac{[4]}{[3]} \qquad
	w = (q-q^{-1})\frac{[2]}{[3]}
\end{equation*}

This recovers the results in \cite{Savage2022} provided we scale by $\lambda = \frac{[4]}{[2]} = q^2+q^{-2}$
and apply the bar involution.

\subsection{Octonion series}
This series is the fourth line in the magic square associated to the octonions.
The distinguishing feature is that $V\cong L$. This is also known as the exceptional
series and is a line in the Vogel plane \cite{Vogel2011}. It has been investigated in
\cite{deligne1996b}, \cite{deligne1996a}, 
\cite{Deligne2002}, \cite[\S~4]{Landsberg2004a}, \cite{Cohen1999}, \cite[Chapter~21]{Cvitanovic2008}, and \cite{Dorey2003}. Quantum dimensions given
in \cite{Westbury2003a} and \cite{Mkrtchyan2017}.

\begin{figure}[ht]
	\begin{tabular}{c|cccccccc}
		$k$ & -4/3 & -1 & -2/3 & 0 & 1 & 2 & 4 & 8 \\ \hline
		& $A_1$ & $A_2$ & $G_2$ & $D_4$ & $F_4$ & $E_6$ & $E_7$ & $E_8$ \\
		& $2\omega_1$ & $\omega_1+\omega_2$ & $\omega_2$ & $\omega_2$ & $\omega_1$ & $\omega_2$ & $\omega_1$ & $\omega_8$
	\end{tabular}
	\caption{Exceptional series}
\end{figure}

\begin{figure}[ht]
	\begin{tabular}{c|ccccc}
		& $V$ & $Y$ & $W$ & $L$ & $X$ \\ \hline
		$A_1$ & $2\omega_1$ & $4\omega_1$ & 0 & $2\omega_1$ & 0 \\
		$A_2\rtimes\fS_2$ & $\omega_1+\omega_2$ & $2\omega_1+2\omega_2$ & $\omega_1+\omega_2$ & $\omega_1+\omega_2$ & $3\omega_1^\ast$ \\
		$G_2$ & $\omega_2$ & $2\omega_2$ & $2\omega_1$ & $\omega_2$ & $3\omega_1$ \\
		$D_4\rtimes\fS_3$ & $\omega_2$ & $2\omega_2$ & $2\omega_1^\ast$ & $\omega_2$ & $\omega_1+\omega_3+\omega_4$ \\
		$F_4$ & $\omega_1$ & $2\omega_1$ & $2\omega_4$ & $\omega_1$ & $\omega_2$ \\
		$E_6$ & $\omega_2$ & $2\omega_2$ & $\omega_1+\omega_6$ & $\omega_2$ & $\omega_4$ \\
		$E_7$ & $\omega_1$ & $2\omega_1$ & $\omega_7$ & $\omega_1$ & $\omega_3$ \\
		$E_8$ & $\omega_8$ & $2\omega_8$ & $\omega_1$ & $\omega_8$ & $\omega_7$ \\
	\end{tabular}
	\caption{Highest weights}\label{series:octonion}
\end{figure}

The equation of the line is $m-n+p=0$ and the parametrisation is
$m=-2$, $n=k+2$, $p=k+4$. The specialisation is
\begin{equation*}
	u \mapsto q^{-2} ,\quad t \mapsto z\,q^{2} ,\quad v \mapsto z\,q^{4}
\end{equation*}

The values of the Casimir are given in Figure~\ref{casimir:octonion}.
\begin{figure}
	\begin{tabular}{ccccc}
		$V$ & $Y$ & $W$ & $L$ & $X$ \\ \hline
		$6k+12$ & $12k+28$ & $10k+16$ & $6k+12$ & $12k+24$
	\end{tabular}
	\caption{Casimir}\label{casimir:octonion}
\end{figure}

The eigenvalues of the braid matrix are given in Figure~\ref{braid:real}.
\begin{figure}
	\begin{tabular}{ccccc}
		$1$ & $L$ & $Y$ & $X$ & $W$ \\ \hline
		$z^{-12}$ & $-q^{-2}z^{-2}$ & $q^2$ & $-1$ & $z^{-6}$
	\end{tabular}
	\caption{Eigenvalues}\label{braid:octonion}
\end{figure}

The dimensions are
\begin{gather*}
	\delta_V = 2\frac{(3k + 7)(5k + 8)}{k+4} \\
	\delta_W = 45 \frac{(2k + 5)(5k + 8)(k + 2)^2}{(k + 4)(k + 6)} \\
	\dim X = 5 \frac{(2k + 5)(3k + 4)(3k + 7)(5k + 8)}{(k + 4)^2} \\
	\dim Y = 90 \frac{(3k + 4)(3k + 7)(k + 2)^2}{(k + 4)^2 (k + 6)} \\
	\dim L = 2 \frac{(3k + 7)(5k + 8)}{k+4}
\end{gather*}
A generalisation is \cite[Theorem~3.2]{Landsberg2002}.

\subsection{Exceptional symmetric spaces}
This is the series in \cite[\S~8.4]{Westbury2015}. Associated to a symmetric space
is a Lie algebra with an involution. The $+1$-eigenspace is a Lie algebra and the
$-1$-eigenspace is an $L$-module, $V$. The symmetric spaces in this series are
characterised by the property that $V$ satisfies the condition \eqref{condition}.

The distinguishing feature is that they all have the same tensor product graph,
\cite{MacKay1991}, \cite{Zhang1991}.
The equation of the line is $n=m$ and it is parametrised by $n=6$, $m=6$, $p=-k$. The Lie
groups in this series are given in Figure~\ref{fig:spinor}.

\begin{figure}[ht]
	\begin{tabular}{c|cccccccc}
$k$ & 2/3 & 1 & 4/3 & 8/3 & 5 & 8 & 10 & 12 \\ \hline
& EVIII & EV & E1 & A1 & & BD1 & FII & EIV \\
& $E_8$ & $E_7$ & $E_6$ & $A_2$ & $OSp(1|2)$ & $D_4$ & $F_4$ & $E_6$ \\
 & $D_8$ & $A_7^\ast$ & $C_4$ & $A_1$ & $A_1$ & $B_3$ & $B_4$ & $F_4$
	\end{tabular}
	\caption{Symmetric spaces}\label{fig:spinor}
\end{figure}

\begin{figure}[ht]
	\begin{tabular}{c|cccccccc}
		$k$ & 5 & 8 & 10 & 12 \\ \hline
		&  & BD1 & FII & EIV \\
		&  $OSp(1|2)$ & $D_4$ & $F_4$ & $E_6$ \\
		 & $A_1$ & $B_3$ & $B_4$ & $F_4$
	\end{tabular}
	\caption{Symmetric spaces}
\end{figure}

The highest weights are given in Figure~\ref{weights:spinor}.
\begin{figure}[ht]
\begin{tabular}{c|ccccc}
& $V$ & $Y$ & $W$ & $L$ & $X$ \\ \hline
$F_4$ & $\omega_4$ & $2\omega_4$ & $\omega_4$ & $\omega_1$ & $\omega_3$ \\
$B_4$ & $\omega_4$ & $2\omega_4$ & $\omega_1$ & $\omega_2$ & $\omega_3$ \\
$B_3$ & $\omega_1$ & $2\omega_4$ & 0 & $\omega_2$ & 0 \\
$A_1$ & $\omega_1$ & 0 & 0 & $2\omega_1$ & 0 \\
$A_1$ & $4\omega_1$ & $8\omega_1$ & $4\omega_1$ & $2\omega_1$ & $6\omega_1$ \\
$C_4$ & $\omega_4$ & $2\omega_4$ & $2\omega_2$ & $2\omega_1$ & $2\omega_3$ \\
$D_8$ & $\omega_8$ & $2\omega_8$ & $\omega_4$ & $\omega_2$ & $\omega_6$ \\
\end{tabular}
\caption{Symmetric spaces}\label{weights:spinor}
\end{figure}

The specialisation is
\begin{equation*}
	u \mapsto q^{-1} ,\quad t \mapsto q^{-1} ,\quad v \mapsto z
\end{equation*}

The values of the Casimir are tentatively given in Figure~\ref{casimir:spinor}.
\begin{figure}
	\begin{tabular}{ccccc}
		$V$ & $Y$ & $W$ & $L$ & $X$ \\ \hline
		$3k-12$ & $6k-20$ & $4k-24$ & $4k-12$ & $6k-24$
	\end{tabular}
	\caption{Casimir}\label{casimir:spinor}
\end{figure}

The eigenvalues of the braid matrix are given in Figure~\ref{braid:real}.
\begin{figure}
	\begin{tabular}{ccccc}
		$1$ & $L$ & $Y$ & $X$ & $W$ \\ \hline
		$z^{6}q^{-12}$ & $-q^{-6}z^{2}$ & $q^{-2}$ & $-1$ & $z^{2}$
	\end{tabular}
	\caption{Eigenvalues}\label{braid:symmetric}
\end{figure}

The dimensions are
\begin{gather*}
	\delta_V = 2 \frac{(10-k)(18-k)}{k} \\
	\delta_W = 9 \frac{(10-k)(12-k)( 15-k)(24-k)}{k^2(6-k)} \\
	\dim X = 2 \frac{(9-k)(10-k)(15-k)(24-k)}{k^2} \\
	\dim Y = 2 \frac{(9-k)(12-k)(15-k)(18-k)}{k(k-6)} \\
	\dim L = (-3) \frac{( 10-k)(12-k)}{k}
\end{gather*}

\subsubsection{Yang-Baxter equation}

Let $V$ be a representation such that $V\otimes V$ is multiplicity free
and the action of the quantised enveloping algebra extends to an
of the affine quantised enveloping algebra. Then the $R$-matrix can be written as
\begin{equation*}
	R(u) = \sum_\mu \rho_\mu(u)\, \pi_\mu
\end{equation*}
The tensor product graph is a method for finding the eigenvalues, $\rho_\mu(u)$,
see \cite{MacKay1991} and \cite{Delius1994}.

The tensor product graph is a bipartite graph. One set of vertices is the set of highest
weights of the symmetric square and the other set of vertices is the set of highest
weights of the exterior square. The two vertices $\mu$ and $\nu$ are connected
if $L(\nu)$ is a composition factor of $L(mu)\otimes L(\theta)$. 
If two vertices $\mu$ and $\nu$ are connected then the eigenvalues satisfy
\begin{equation*}
	(uq^C-u^{-1}q^{-C})\rho_\mu(u) = (uq^C-u^{-1}q^{-C})\rho_\nu(u)
\end{equation*}

For this series the tensor product graph is
\begin{equation*}
	\begin{tikzcd}
	I \arrow[r, dash] & L \arrow[r, dash] & Y \arrow[r, dash] & X \arrow[r, dash] & W
\end{tikzcd}
\end{equation*}

\begin{defn}\label{defn:Rmatrix} The $R$-matrix is
	\begin{gather*}
		 [x-2k+3][x-k+2][x+1][x+k] \\
		+ [-x-2k+3][x-k+2][x+1][x+k]  \, \pi_L \\
		+ [-x-2k+3][-x-k+2][x+1][x+k]   \, \pi_Y \\
		+ [-x-2k+3][-x-k+2][-x+1][x+k]   \, \pi_X \\
		+ [-x-2k+3][-x-k+2][-x+1][-x+k]   \, \pi_W
	\end{gather*}
\end{defn}

\begin{lemma}
\begin{equation*}
	R(u)\,R(u^{-1}) =  [x-k+3][x+3][x-k/2][x-1][-x-k+3][-x+3][-x-k/2][-x-1]
\end{equation*}
\end{lemma}
\begin{proof}
This follows from Definition~\ref{defn:Rmatrix}.
\end{proof}

\begin{prop}\label{prop:YB} The $R$-matrix, $R(u)\in \End([2])[u^{\pm 1}]$ is given by
	\begin{gather*}
		R(u) = [x+2k-3][x+k][x+k-3]\left(\frac{uS^{-1}-u^{-1}S}{q-q^{-1}}\right) \\
		+ \frac{[k-2]}{[k-3]}[x+k][x+k-3]\;U \\
		+ \frac{[k-4][k]}{[k-3]}[x+2k-3][x+k-3]\;K \\
		-[k-4][x+2k-3][x+k]\;H
	\end{gather*}
\end{prop}

\begin{proof}
	This is a direct calculation using computer algebra.
\end{proof}
The leading coefficient is the coefficient of $u^4$ and this is $\frac{z^{-2}q^{6}}{(q-q^{-1})^4}S$.
Similarly, the trailing coefficient is the coefficient of $u^{-4}$ and this is $\frac{z^{2}q^{-6}}{(q-q^{-1})^4}S$.

We can also regard $R(u)$ as a function from invertible scalars to $\End([2])$ by substituting for $u$.
Then we have:
Note that 
\begin{align*}
	R(1) & =  [2k-3][k-2][k] \\
	R(z^{-2}q^3) &= [k-3][3k-5][2k-4][4k-6]\;U \\
	R(z^{-1}) &= -[2][k-3][k+1][2k]\;K \\
	R(z^{-1}q^3) &=-[2][k-3][k+1][2k]\;H
\end{align*}

\begin{lemma}[Unitarity] $R(1) =[2k-3][k-2][k]$.
\end{lemma}
\begin{proof}
	Substituting in the definition gives
	\begin{equation*}
		R(1) = [9][6]\left([3]\left(\frac{S-S^{-1}}{q-q^{-1}}\right)
		+ [4]\;U
		- [4]\;K
		+ [4]\;H	\right)
	\end{equation*}
	Hence it is sufficient to show
	\begin{equation*}
		[4] = [3]\left(\frac{S-S^{-1}}{q-q^{-1}}
		+ [4]\;U
		- [4]\;K
		+ [4]\;H	\right)
	\end{equation*}
	This is the skein relation.
\end{proof}

\begin{prop}[Crossing symmetry] The $R$-matrix satisfies
\begin{equation*}
	\Rot(R(u)) = R(u^{-1}z^{-2}q^{6})
\end{equation*}
\end{prop}
\begin{proof}
This is a direct calculation using computer algebra.
\end{proof}

\begin{defn}
	Define $R_1(u), R_2(u)\in \End([3])[u,u^{-1}]$ by $R_1(u) = R(u)\otimes 1$ and $R_2(u)=1\otimes R(u)$.
\end{defn}

\begin{prop}[Yang-Baxter] The $R$-matrix, $R(u)\in \End([2])[u,u^{-1}]$ satisfies
	\begin{equation*}
		R_1(u)R_2(uv)R_1(v) = R_2(v)R_1(uv)R_2(u)
	\end{equation*}
	This relation holds in $\End([3])[u,v,(uv)^{-1}]$
\end{prop}
We require this relation to be satisfied since it holds in $\End(\otimes^3V)[u,v,(uv)^{-1}]$
after applying the evaluation functor.

\subsection{Orthogonal series}
This is a degenerate line. Taking the specialisation and quotienting by the negligible ideal
gives the Birman-Wenzl category. The distinguishing feature is that $W=0$ and $X=0$.

The equation of the line is $2m+2n+p=0$ and the parametrisation is
$m=3$, $n=N-4$, $p=2-2N$.
\begin{figure}[ht]
	\begin{tabular}{c|ccccc}
		& $V$ & $Y$ & $W$ & $L$ & $X$ \\ \hline
		$B_n$ & $\omega_1$ & $2\omega_1$ & 0 & $\omega_2$ & 0 \\
		$D_n$ & $\omega_1$ & $2\omega_1$ & 0 & $\omega_2$ & 0 \\
	\end{tabular}
	\caption{Orthogonal}\label{series:orthogonal}
\end{figure}

\subsection{$SL(2)$}
This is line $2m+3n+2p=0$. This is invariant under reflection.
The specialisation is
\begin{equation*}
	t\mapsto q^{-2},\quad u\mapsto z^3,\quad v\mapsto z^{-3}q^3
\end{equation*}

All points on this line give the quantised enveloping algebra of $SL(2)$
with quantum parameter $q^4$. The representation $V$ is the two dimensional
irreducible representation taken as an odd super vector space.


\bibliographystyle{alphaurl}
\bibliography{common}

\end{document}